\documentclass{article}[12pt]
\usepackage{amsmath,amsfonts,amsthm}
\usepackage{amssymb}

\usepackage{epsfig}

\usepackage{subfigure}
\usepackage{graphics}
\usepackage{graphicx}
\usepackage{color}
\usepackage{soul}
\usepackage[makeroom]{cancel}

\usepackage{latexsym}
\usepackage{graphics}
\usepackage{amsmath}
\usepackage{amsthm}
\usepackage{xspace}
\usepackage{amssymb}
\usepackage{color}
\usepackage{cite}
\usepackage{latexsym}
\usepackage{graphics}
\usepackage{amsmath}
\usepackage{amsthm}
\usepackage{xspace}
\usepackage{amssymb}
\usepackage{epsfig}
\usepackage{comment}

\newcommand{\beql}[1]{\begin{equation}\label{#1}}
\newcommand{\eeql}{\end{equation}}
\newcommand{\eqn}[1]{(\ref{#1})}

\newcommand{\yuan}[1]{#1}

\newcommand{\R}{\mathbb{R}}
\newcommand{\pr}{\mathbb{P}}

\newcommand{\ci}{{\cal I}}

\newcommand{\ck}{{\cal K}}
\newcommand{\cx}{{\cal X}}
\newcommand{\cm}{{\cal M}}

\newcommand{\ch}{{\cal H}}

\newcommand{\bk}{\boldsymbol{k}}
\newcommand{\bx}{\boldsymbol{x}}
\newcommand{\bX}{\boldsymbol{X}}
\newcommand{\bbeta}{\boldsymbol{\eta}}
\newcommand{\bnu}{\boldsymbol{\nu}}
\newcommand{\be}{\boldsymbol{e}}
\newcommand{\bu}{\boldsymbol{u}}
\newcommand{\bZero}{\boldsymbol{0}}

\newtheorem{thm}{Theorem}
\newtheorem{lem}[thm]{Lemma}

\newtheorem{definition}[thm]{Definition}
\newtheorem{conjecture}[thm]{Conjecture}

\begin{document}

\title{Asymptotic optimality of a greedy randomized algorithm in a large-scale
service system with general packing constraints
}

\author
{
Alexander L. Stolyar \\
Bell Labs, Alcatel-Lucent\\
600 Mountain Ave., 2C-322\\
Murray Hill, NJ 07974 \\
\texttt{stolyar@research.bell-labs.com}
\and
Yuan Zhong\\
Columbia University\\ 
500 W. 120 St, Mudd 344\\
New York, NY 10027\\
\texttt{yz2561@columbia.edu}
}

\date{\today}

\maketitle

\begin{abstract}

We consider a service system model primarily
motivated by the problem of efficient assignment of virtual machines to physical host machines in 
a network cloud, so that the number of occupied hosts is minimized.

There are multiple  types of  arriving 
customers,  where a customer's mean service time  depends 
on its type. 
There is an infinite number of servers. 
Multiple customers can be placed for service into 
one server, subject to general ``packing'' constraints.
Service times of different customers are independent, even if served simultaneously
by the same server.  
Each new arriving customer is placed for service immediately, either
 into a server already serving other customers 
(as long as packing constraints are not violated)
or into an idle server. 
After a service completion, each customer leaves its server and the system.

We propose an extremely simple and easily implementable
customer placement algorithm, called {\em Greedy-Random} (GRAND). 
It places each arriving customer uniformly at random into either one of the already occupied servers (subject to packing constraints)
or one of the so-called {\em zero-servers}, which are empty servers designated to be available to new arrivals.
One instance of GRAND, called GRAND($aZ$), where $a\ge 0$ is a parameter, is such that the number of zero-servers
at any given time $t$ is $aZ(t)$, where $Z(t)$ is the current total number of customers in the system.
We prove that GRAND($aZ$) with $a>0$ 
is asymptotically optimal, as the  customer arrival rates grow to infinity and $a\to 0$,
in the sense of minimizing the total number of occupied servers in steady state. 
In addition, we study by simulations various versions of GRAND 
and observe the dependence of convergence speed and steady-state performance 
on the number of zero-servers.

\end{abstract}

\newpage

\section{Introduction}
\label{sec-intro}

The model in this paper is primarily motivated by the 
problem of efficient real-time assignment (``packing'') of virtual machines 
(VM)
to physical host machines in a cloud  computing facility   \cite{Gulati2012},
so that the number of occupied hosts is minimized. 
We refer readers to \cite{Gulati2012}
for a general discussion of resource allocation problems in cloud computing,
including virtual machine assignment.

There is an infinite number of servers,  which   serve customers belonging to one of a finite number of types. 
Multiple customers can be placed for service into 
one server, subject to general ``packing'' constraints.
Service times of different customers are independent, even if served simultaneously
by the same server. 
Each new arriving customer is placed for service immediately, either
 into a server already serving other customers 
(as long as packing constraints are not violated)
or into an idle server. 
After service completion, each customer leaves its server and the system --
this is the key feature that separates this model from classical 
bin packing models (see, e.g., \cite{Csirik2006,Bansal2009}).

We are interested in finding algorithms for real-time assignment
of arriving customers to servers, so that the total number 
of  occupied servers in the stationary regime is minimized.
For the purpose  of practical implementation,
it is very desirable that such an algorithm be  as simple as possible,
and require  as little information about the system parameters
and state as possible.

Papers \cite{St2012,StZh2012} introduce and study
different {\em Greedy} algorithms,
and prove their asymptotic optimality, as the customer  arrival rates grow to 
infinity. A Greedy algorithm does not need to know the customer  arrival rates
or service times, and makes placement decisions based on the current
system state only. It does, however, need to keep track of the numbers $X_{\bk}$
of servers being in different packing configurations $\bk$, where a configuration $\bk$ describes  the numbers of   customers of each type that a server contains. 
(The algorithm in \cite{St2012} greedily pursues the minimization of the convex objective $\sum_{\bk} X_{\bk}^{1+\alpha}$,
with $\alpha>0$ being a parameter of the algorithm; this objective approximates the total number of occupied servers
$\sum_{\bk} X_{\bk}$ when $\alpha$ is small. In \cite{StZh2012} a different convex modification of $\sum_{\bk} X_{\bk}$ is used.)
A practical implementation of the {\em Greedy} algorithms in \cite{St2012,StZh2012} may be difficult,
if the number of different configurations $\bk$ is large. 
 In many practical scenarios, this number is indeed very large.

In this paper,  
we propose an extremely simple and easily implementable
customer placement algorithm, called {\em Greedy-Random} (GRAND).  
In particular, we study one instance of this algorithm, called GRAND($aZ$),
where $a\ge 0$ is a parameter. Our main results prove that  GRAND($aZ$) with $a>0$
is asymptotically optimal, as the customer  arrival rates grow to infinity
and $a\to 0$. 

Just like previous Greedy algorithms referred to above, GRAND does not need to know 
the customer arrival rates or service times. {\em Moreover, GRAND($aZ$) with $a>0$ does {\em not}
need to keep track of the exact configurations of the servers!} 
For each $a > 0$, GRAND($aZ$) works basically as follows:\\
(a) it keeps track of which servers are non-empty (but not of their exact configurations);\\
(b) it keeps track of the total number $Z$ of customers 
in the system (which is usually not difficult to track);\\
(c) at any time, it
designates the subset of $aZ$ of empty servers as ``zero-servers'' (recall that there is always an infinite number of available empty servers);\\
(d) upon arrival of a customer of {\em any type}, it chooses  a server  uniformly 
 at random  among all  occupied   and zero-servers, and places the customer
into this server, unless this will violate the packing constraints;\\
(e) step (d) is repeated until the placement succeeds, which occurs
in finite time. 

We see that GRAND($aZ$)  is as simple as any algorithm can posibly be, in terms of the amount
of system state information it needs to track.
It is perhaps counterintuitive that an algorithm as simple and   ``blind''
as GRAND($aZ$) can result in an asymptotically optimal packing in steady state.
Indeed, the intuition from the previous work on classical bin packing (where customers never leave the system)
suggests that ``naive greedy packing algorithms cannot possibly work well," except in very special cases. 
Our results show that this intuition is false, at least in principle, when we consider systems with customer departures.
In a system with departures, customer placement decisions are in a sense ``reversible," because each customer leaves eventually,
and therefore ``bad" placement decisions can be ``undone" over time. 
Still, the fact that even a procedure as ``naive" as GRAND($aZ$) can be optimal, is surprising.
We consider this insight to be the main contribution of this paper.

In addition to the theoretical results on the GRAND($aZ$) optimality,
we study different versions of GRAND via simulations. 
We observe a tradeoff between steady-state performance and 
convergence speed, 
as we vary the number of zero-servers. 
A major conclusion from these experiments is that
for systems of reasonable (not extremely large) scale,
the GRAND($0$) algorithm, i.e. GRAND($aZ$) with $a=0$, 
has the best performance. This makes GRAND($0$) 
very attractive for practical applications.

\subsection{Brief review of previous work}

Our work is related to several lines of previous research.
First, there is an extensive literature on classical 
{\em online stochastic bin packing} problems (see e.g. \cite{Csirik2006,Bansal2009} for good recent reviews). In such problems,
random-sized items arrive   to the system and need to be placed 
according to an online algorithm
 into finite-size bins.   \yuan{Arriving} items never leave or move between bins,
and the typical objective is to minimize the number of occupied bins.
A very recent paper \cite{GR2012}
provides new strong results for this classical setting. It also contains
some heuristics and simulations for the case with item departures,
which is a special case of our model.
The term {\em vector-packing} (cf. \cite{St2012}) refers to multi-dimensional bin packing, 
when bins and item sizes are vectors. We remark that the packing constraints in our model
include vector-packing constraints as a special case.

Another line of previous work is on 
bin packing {\em service} systems (see e.g. \cite{Gamarnik2004}
and references therein). In such systems,  random-sized 
items (or customers) \yuan{arrive according to a stochastic process.}  
 \yuan{Arriving customers} need to be served by a bin (\yuan{or} server), and leave after a random service time. 
A server can simultaneously process multiple customers as long as 
they can simultaneously fit into it. Customers waiting for service are queued
\yuan{, and} a typical problem is to determine
the maximum throughput under a given packing algorithm
for assigning  customers for service. 
 Our model is similar to these systems 
 \yuan{in the sense of considering customer departures},
 \yuan{but in our model,} there is an infinite
number of available bins (\yuan{or} servers), which means there are no queues or stochastic stability
issues.

A recent line of work on bin packing service systems with multiple servers,
which is also motivated by real-time VM allocation problems,
includes \cite{Jiang2012, Maguluri2012, Maguluri2013, GSW2012}.
In \cite{Jiang2012}, authors address a real-time VM allocation problem, 
which in particular includes packing constraints.
 \yuan{The} approach of \cite{Jiang2012} is close in spirit to Markov Chain algorithms used in 
combinatorial optimization.
Papers \cite{Maguluri2012,Maguluri2013} are
 concerned mostly with online algorithms   that maximize 
 throughput of a queueing system 
with a finite number of bins, based on current queue lengths. 
In \cite{GSW2012}, the approach is to dynamically solve the underlying
optimization problem by employing a shadow queueing system, 
and guide actual packing by the current solution.

The model in this paper is \yuan{the} same as that studied in \cite{St2012, StZh2012},
but the algorithm that we consider here, GRAND, is   \yuan{drastically} different from the versions
of Greedy algorithms in \cite{St2012, StZh2012} -- we already discussed this
in Section~\ref{sec-intro}.

The asymptotic regime in this paper is such that   arrival rates 
and the average number of occupied servers
scale up to infinity.
In this respect, our work is related to the extensive literature on queueing systems in the 
{\em many servers} regime. (See e.g. \cite{ST2010_04} for an overview.)
However, packing constraints 
are not present in work\yuan{s} 
prior to \cite{St2012, StZh2012}
on the many servers regime, to the best of our 
knowledge.

Finally, we mention the works \cite{zachary2000,kelly1986,kelly1991}, which study the large-scale asymptotics
of an uncontrolled loss system. In this paper, as one of the tools, we use an ``entropy-like'' objective/Lyapunov function, which is similar in form
to that employed in the above papers. We note, however, that there is {\em no} close analogy between our main results and those of \cite{zachary2000,kelly1986,kelly1991}, due to substantial differences in the models and algorithms.
(See Remark 1 in Section~\ref{sec-model}.)

\subsection{Layout of the rest of the paper}

We formally define the model, the GRAND algorithm, and state the main
asymptotic optimality results in Section~\ref{sec-model}.
Sections~\ref{sec-opt-conv}  \yuan{and} \ref{sec-fsp-dynamics} contain the proofs,
with the central part being the analysis of convergence of fluid sample paths
in Section~\ref{sec-fsp-dynamics}. A probabilistic interpretation
of the \yuan{steady-state solution of} GRAND(aZ)  is given in Section~\ref{sec-opt-point-interpretation}.
In Section~\ref{sec-sim} we present and discuss simulation experiments.
Concluding remarks and discussion,  including some conjectures, are in Section~\ref{sec-further-work}.

\subsection{Basic notation and conventions}
\label{subsec-notation}

Sets of real and real non-negative numbers are denoted by $\R$ and $\R_+$, respectively.
The standard Euclidean norm of a vector $\bx\in \R^n$ is denoted by $\|\bx\|$. The distance from vector $\bx$ to a set $U$ in a Euclidean space is denoted by 
$d(\bx,U)=\inf_{\bu\in U} \|\bx-\bu\|$. 
$\bx \to \bu \in \R^n$ means ordinary convergence in $\R^n$,
and $\bx \to U \subseteq \R^n$
means $d(\bx,U) \to 0$. $\be_i$ is the $i$-th coordinate unit vector in $\R^n$.
Symbol $\implies$
denotes convergence in distribution of random variables taking values in space $\R^n$
equipped with the Borel $\sigma$-algebra. The
abbreviation {\em w.p.1} means convergence {\em with probability 1}.
We often write $x(\cdot)$ to mean the function (or random process) $\{x(t),~t\ge 0\}$.
Abbreviation {\em u.o.c.} means 
{\em uniform on compact sets} convergence of functions, with the argument (usually in 
$[0,\infty)$) determined by the context.
We often  write $\{x_{\bk}\}$ to mean the vector $\{x_{\bk}, ~\bk\in\ck\}$, with the set of indices
$\ck$ determined by the context.
The cardinality of a finite set $\mathcal{N}$ is 
$|\mathcal{N}|$. Indicator function $I\{A\}$ for a condition $A$ is equal to $1$
if $A$ holds and $0$ otherwise.  $\lceil \xi \rceil$ denotes the smallest integer
\yuan{greater than or equal to} $\xi$ 
\yuan{, and $\lfloor \xi \rfloor$ denotes the largest integer 
smaller than or equal to $\xi$.}

For a finite set of scalar functions $f_n(t), ~t\ge 0$, $n\in\mathcal{N}$, a point $t$ is called
{\em regular} if for any subset $\mathcal{N}' \subseteq \mathcal{N}$ the 
derivatives
$$
\frac{d}{dt} \max_{n\in\mathcal{N}'} f_n(t) ~~\mbox{and}~~ 
\frac{d}{dt} \min_{n\in\mathcal{N}'} f_n(t)
$$
exist.   (To be precise, we require that each derivative is proper: both left and right derivatives exist and are equal.)
 In this paper, we use bold letters to represent vectors, 
and plain letters to represent scalars.

\section{Model and main results}
\label{sec-model}

The model we study in this paper is the same as in \cite{St2012,StZh2012}.
We consider a service system with $I$  types 
\yuan{of arriving customers}, indexed by $i \in \{1,2,\ldots,I\} \equiv \ci$. 
\yuan{Customers of type $i$ arrive as an independent Poisson process of} rate $\Lambda_i >0$. The service time of a type-$i$ customer is an exponentially distributed random variable with mean $1/\mu_i$.
All  \yuan{customer inter-arrival times} and  service times are mutually independent.
There is an infinite number of \yuan{homogeneous} servers. Each server can potentially serve more than one customer simultaneously, subject to the following general packing constraints. We say that a vector $\bk = \{k_i, ~i\in \ci\}$ with non-negative integer components
is a server {\em configuration}, if a server can simultaneously serve a combination of customers \yuan{of different types} given by \yuan{the} vector $\bk$. The set of all configurations is finite, and is denoted by $\bar\ck$.
We assume that $\bar\ck$ satisfies the {\em monotonicity} condition: if 
$\bk\in \bar\ck$, then all ``smaller'' configurations $\bk'$ \yuan{with} $\bk'\le \bk$ \yuan{component-wise}, belong to $\bar\ck $ as well. Without loss of generality, assume that $\be_i \in \bar\ck$ for all types $i$, where $\be_i$ is the $i$-th coordinate unit vector (otherwise, type-$i$ customers cannot be served at all).
By convention, the component-wise zero vector $\bk=\bZero$ belongs to $\bar\ck$ -- this is the configuration of an empty server. We denote by $\ck=\bar\ck \setminus \{\bZero\}$ 
the set of configurations {\em not} including the zero configuration. 

An important feature of the model is that simultaneous service does {\em not} affect the service time
distributions of individual customers. In other words, the service time of a customer is unaffected by whether or not there are other customers served simultaneously by the same server. Each arriving customer is immediately placed for service in one of the servers; it can be ``added'' to an empty or occupied server, as long as the configuration feasibility constraint is not violated,
i.e. it can be added to any server whose \yuan{current} configuration $\bk\in\bar\ck$  is such that $\bk+\be_i \in \ck$. 
When the service of a type-$i$ customer by the server in configuration $\bk$ is completed,
the customer leaves the system and the server's configuration changes to $\bk-\be_i$.
Denote by $X_{\bk}$ the number of servers with configuration $\bk\in \ck$. The system state is then the vector $\bX = \{X_{\bk}, ~\bk\in \ck\}$. 

A {\em placement algorithm} (or packing rule) determines where an arriving customer is placed, as a function of the current system state $\bX$. Under any well-defined placement algorithm, 
\yuan{the process $\{\bX(t), t\ge 0\}$} is a continuous-time Markov chain 
\yuan{with a countable state space}. It is easily seen to be irreducible and positive recurrent: the positive recurrence follows from the fact
that the total number $Y_i(t)$ of type-$i$ customers in the system is independent from the  \yuan{placement algorithm}, and its stationary distribution is 
Poisson with mean $\Lambda_i/\mu_i$; we denote by $Y_i(\infty)$ the random value of $Y_i(t)$ in steady-state -- it is, therefore, a Poisson random 
variable with mean $\Lambda_i/\mu_i$.
Consequently, the process $\{\bX(t), ~t\ge 0\}$ has a unique stationary distribution; let $\bX(\infty) = \{X_{\bk}(\infty), \bk \in \ck \}$ be the random system state $\bX(t)$ in stationary regime.

We are interested in finding a 
\yuan{placement algorithm that minimizes} 
the total number of  \yuan{occupied} servers $\sum_{\bk\in\ck} X_{\bk}(\infty)$
in the stationary regime. 


\begin{definition}[Greedy-Random (GRAND) algorithm]\label{df:grand}
At any given time $t$, there is a designated finite set of $X_{\bZero}(t) \ge 0$ empty servers, 
called {\em zero-servers}.  
(Recall that there is always an infinite number of available empty servers.)
$X_{\bZero}(t) = f(\bX(t))$ is a given fixed function 
of the state $\bX(t)$. (Recall that $\bX(t)= \{X_{\bk}(t), ~\bk\in \ck\}$ describes the quantities of servers in different non-zero configurations.) \\
1. A new customer, 
say of type $i$, arriving at time $t$
is placed into a server chosen
 uniformly \yuan{at random} among the zero-servers and those occupied servers,
where it can still fit. In other words,
the total number of servers available to a type-$i$ arrival at time $t$ is
$$
X_{(i)}(t) \doteq X_{\bZero}(t) + \sum_{\bk\in \ck:~\bk+\be_i\in \ck} X_{\bk}(t).
$$
If $X_{(i)}(t)=0$, the customer is placed into an empty server.\\
2. It is important
that immediately after any arrival or departure at time $t$, 
the value of $X_{\bZero}(t+0)$ 
is reset (if necessary) to  $f(\bX(t+0))$.
\end{definition}

We now define a special version of GRAND, which will be the main focus of this paper.

\begin{definition}[GRAND($aZ$)  algorithm]\label{df:grand-az}
A special case of the GRAND algorithm, with 
the number of zero-servers depending only on the 
total number of customers $Z=\sum_i \sum_{\bk} k_i X_{\bk}$  as
$$
X_{\bZero}=\lceil a Z \rceil,
$$ 
where $a\ge 0$ is a parameter, will be called GRAND($aZ$) algorithm.
\end{definition}

We emphasize that a GRAND algorithm is extremely 
simple and easy to implement, as long as $X_{\bZero}$ is a relatively simple
 function of $\bX$. For example, the following procedure implements  GRAND($aZ$) with $a>0$.
Suppose \yuan{that}, upon a customer arrival, $Z>0$. (The case
$Z=0$ is trivial -- the arriving customer is 
placed into an empty server.)
We pick  uniformly \yuan{at random} any server among 
\yuan{all} zero-servers and occupied servers.
If the customer fits into the  \yuan{chosen}
server, place it there; if not, repeat the uniform random choice until
an available server is found. We remark that the procedure stops eventually, 
because $a>0$ and $Z>0$ guarantees
$X_{\bZero}\ge 1$, and then $X_{(i)}\ge 1$.
{\em There is no need to keep track of the exact configurations of   individual
servers.  Only the set of occupied servers and the value of $Z(t)$ need to be known.}

We now define the asymptotic regime. Let $r\to\infty$ be a positive scaling parameter.
 More specifically, assume that $r\ge 1$, and $r$ increases to infinity along a discrete
sequence. 
Customer arrival rates scale linearly with $r$; namely,
for each $r$, $\Lambda_i = \lambda_i r$, where $\lambda_i$ are fixed positive parameters.
Let $\bX^r(\cdot)$ be the process associated with \yuan{a} system with parameter $r$, and 
\yuan{let} $\bX^r(\infty)$ be the (random) system state in the stationary regime.
For each $i$, denote by $Y^r_i(t) \equiv \sum_{\bk\in\ck} k_i X^r_{\bk}(t)$ the total number
of customers of type $i$. Since arriving customers are placed for service immediately
and their service times are independent \yuan{of each other and} of the rest of the system, 
 $Y^r_i(\infty)$ is a Poisson random variable with mean $r \rho_i$, where $\rho_i\equiv \lambda_i/\mu_i$.
Moreover, $Y^r_i(\infty)$ are independent across $i$.
Since the total number of occupied servers is no greater than the total number 
of 
customers, 
$\sum_{\bk\in\ck} X_{\bk}^r(t) \le Z^r(t)\equiv \sum_i Y^r_i(t)$, we have a simple upper bound 
on the total number of occupied servers in steady state,
$\sum_{\bk\in\ck}  X_{\bk}^r(\infty) \le Z^r(\infty) \equiv \sum_i Y^r_i(\infty)$,
where $Z^r(\infty)$ is a Poisson random variable with mean $r \sum_i \rho_i$.
Without loss of generality, from now on  we assume $\sum_i \rho_i=1$.
This is equivalent to rechoosing the parameter $r$ to be $r \sum_i \rho_i$.

By convention, we will {\em not} include $X_{\bZero}^r(t)$ into $\bX^r(t)=\{X_{\bk}^r(t),~\bk\in \ck\}$.

From this point on in this section, as well as in Sections~\ref{sec-opt-conv} --
\ref{sec-opt-point-interpretation},
we consider only the  GRAND($aZ$) \yuan{algorithms} with $a>0$.

The {\em fluid-scaled} process is $\bx^r(t)=\bX^r(t)/r$, $t \in [0, \infty)$. 
We also define $\bx^r(\infty) = \bX^r(\infty)/r$.
For any $r$, $\bx^r(t)$ takes values in the  non-negative orthant $\R_+^{|\ck|}$.
Similarly, $y^r_i(t)=Y^r_i(t)/r$, $z^r(t)=Z^r(t)/r$, $x^r_{\bZero}(t)=X^r_{\bZero}(t)/r$ and 
$x^r_{(i)}(t)=X^r_{(i)}(t)/r$,  for  $t \ge 0$ and $t=\infty$.
Since $\sum_{\bk\ne 0}  x_{\bk}^r(\infty) \le z^r(\infty)=Z^r(\infty)/r$,
we see that the random variables
$(\sum_{\bk\ne 0} x_{\bk}^r(\infty))$ are uniformly integrable in $r$.
This in particular implies that the sequence of distributions of $\bx^r(\infty)$ is tight,
and therefore there always exists a limit \yuan{$\bx(\infty)$} in distribution\yuan{, so that} $\bx^r(\infty)\implies \bx(\infty)$,
along a subsequence of $r$. 

The limit (random) vector $\bx(\infty)$ satisfies the following conservation laws:
\beql{eq-cons-laws}
\sum_{\bk\in\ck} k_i x_{\bk}(\infty) \equiv y_i(\infty) = \rho_i, ~~\forall i,
\end{equation}
\yuan{which,} in particular, \yuan{implies that}
\beql{eq-cons-laws2}
z_i(\infty)\equiv \sum_i y_i(\infty) \equiv \sum_i \rho_i
= 1.
\end{equation}
Therefore, the values of $\bx(\infty)$ are confined to the convex compact 
$(|\ck|-I)$-dimensional 
polyhedron
$$
\cx \equiv \{\bx\in \R_+^{|\ck|} ~|~  \sum_{\bk} k_i x_{\bk} = \rho_i, ~\forall i\in\ci \}.
$$
We will slightly abuse notation by using symbol $\bx$ for a generic element of $\cx$;
while $\bx(\infty)$, and later $\bx(t)$, refer to random  \yuan{vectors} taking values in $\cx$.

Also note that under GRAND($aZ$), $x^r_{\bZero}(\infty) \implies x_{\bZero}(\infty)=a z(\infty)=a$, as $r\rightarrow \infty$.

The asymptotic regime and the associated basic properties \eqn{eq-cons-laws}
and \eqn{eq-cons-laws2} hold {\em for any  \yuan{placement algorithm}}. 
 Indeed, \eqn{eq-cons-laws}
and \eqn{eq-cons-laws2} only depend on the already mentioned  fact 
that all $Y_i^r(\infty)$ are mutually independent
Poisson  random variables
with means $\rho_i r$.

Consider the following problem of minimizing the number of occupied servers,
on the fluid scale:
$\min_{\bx\in \cx} \sum_{\bk\in\ck} x_{\bk}$. It is a
linear program:
\beql{eq-opt}
\min_{\bx\in\R_+^{|\ck|}} \sum_{\bk\in\ck} x_{\bk} ,
\end{equation}
subject to
\beql{eq-cons-laws222}
\sum_{\bk\in\ck} k_i x_{\bk} = \rho_i, ~~\forall i.
\end{equation}
Denote by $\cx^* \subseteq \cx$ the set of its optimal solutions.

For future reference, we record the following observations and notation. 
Using the monotonicity of $\bar{\ck}$, it is easy to check that if in the LP 
\eqn{eq-opt}-\eqn{eq-cons-laws222} we replace equality constraints
\eqn{eq-cons-laws222} with the inequality constraints
\beql{eq-cons-laws222ge}
\sum_{\bk\in\ck} k_i x_{\bk} \ge \rho_i, ~~\forall i,
\end{equation}
the new LP \eqn{eq-opt},\eqn{eq-cons-laws222ge} has same optimal value, and its set of the optimal solutions 
$\cx^{**}$ contains $\cx^*$, or more precisely, $\cx^* = \cx^{**} \cap \cx$.
From here,  using Kuhn-Tucker theorem,
$\bx \in \cx^*$ if and only if there exists a vector $\bbeta=\{\eta_i, ~i\in \ci\}$
of Lagrange multipliers, corresponding to the
inequality  constraints \eqn{eq-cons-laws222ge}, such that the following conditions hold:
\beql{eq-dual-1}
\bx\in \cx,
\end{equation}
\beql{eq-dual-111}
\eta_i \ge 0, ~~\forall i \in \ci,
\end{equation}
\beql{eq-dual-2}
\sum_i k_i \eta_i \le 1, ~~\forall \bk \in \ck,
\end{equation}
\beql{eq-dual-3}
\mbox{for any $\bk \in \ck$, ~~ condition} ~\sum_i k_i \eta_i < 1 ~\mbox{implies}~x_{\bk}=0.
\end{equation}
Vectors $\bbeta$  satisfying \eqn{eq-dual-1}-\eqn{eq-dual-3} for some $\bx \in \cx$
 are optimal solutions to the problem dual to LP \eqn{eq-opt},\eqn{eq-cons-laws222ge}.
 They form a convex set, which we denote by $\ch^*$; it is easy to check that $\ch^*$ is compact.

For each $a>0$, denote 
\beql{eq-L-def}
L^{(a)}(\bx) = [-\log a]^{-1}\sum_{\bk\in\ck} x_{\bk} \log [x_{\bk} c_{\bk} /(e a)],
\end{equation}
where $c_{\bk} \doteq \prod_i k_i !$, $0!=1$. 
Throughout the paper we will use notation $b= -\log a$.
We then have 
\beql{eq-L-partial}
(\partial / \partial x_{\bk}) L^{(a)}(\bx) =(1/b)  \log [x_{\bk} c_{\bk} / a ], ~~\bk\in \ck.
\end{equation}
Note that if we adopt a convention that
\beql{eq-formal-deriv-zero}
(\partial / \partial x_{\bZero}) L^{(a)}(\bx)|_{x_{\bZero}=a} = 0,
\end{equation}
then \eqn{eq-L-partial} is valid for $\bk=\bZero$ and $x_{\bZero}=a$, which will be useful later.

The function $L^{(a)}(\bx)$ is strictly convex in $\bx\in\R_+^{|\ck|}$.
Consider the problem $\min_{\bx\in \cx} L^{(a)}(\bx)$. It is the
 following convex optimization problem:
\beql{eq-opt-grand}
\min_{\bx\in\R_+^{|\ck|}} L^{(a)}(\bx),
\end{equation}
subject to
\beql{eq-cons-laws-grand}
\sum_{\bk\in\ck} k_i x_{\bk} = \rho_i, ~~\forall i.
\end{equation}
Denote by $\bx^{*,a} \in \cx$ its unique optimal solution. 
Using \eqn{eq-L-partial} it is easy to check that $x^{*,a}_{\bk}>0$ for all ${\bk} \in \ck$.
There exists a vector $\bnu^{*,a} = \{\nu_i^{*,a}, ~ i\in \ci\}$ of  Lagrange multipliers for 
the constraints \eqn{eq-cons-laws-grand}, such that $\bx^{*,a}$ solves problem
$$
\min_{\bx\in \R_+^{|\ck|}} L^{(a)}(\bx) + \sum_i \nu_i^{*,a} (\rho_i - \sum_{\bk\in\ck} k_i x_{\bk}).
$$
By setting partial derivatives of the objective function to zero, we see that $(1/b)  \log [x_{\bk}^{*,a} c_{\bk} / a ] - \sum_i \nu_i^{*,a} k_i = 0$, $ \bk \in \ck$. 
Therefore, $\bx^{*,a}$ has the product form
\beql{eq-grand-product}
x^{*,a}_{\bk} = \frac{a}{c_{\bk}} \exp \left[b \sum_i k_i \nu_i^{*,a}\right] = \frac{1}{c_{\bk}} a^{1-\sum_i k_i \nu_i^{*,a}} , ~~ \bk \in \ck.
\end{equation}
This in particular implies that the Lagrange multipliers $\nu_i^{*,a}$ are unique 
and are equal to $1 - \log (x_{\be_i}^{*,a})/\log a$, by considering \eqref{eq-grand-product} for $\be_i$, $i \in \ci$; note also that they
can have any sign (not necessarily non-negative). We can summarize: {\em
a point $\bx\in \cx$ is the optimal solution to \eqn{eq-opt-grand}-\eqn{eq-cons-laws-grand}
(that is $\bx=\bx^{*,a}$) if and only if it has a product form representation
\eqn{eq-grand-product} for some vector $\bnu^{*,a}$.} (The 'only if' part we just proved, and the 'if' follows
from Kuhn-Tucker theorem.)

Our {\bf main results} are as follows. Speaking informally,
they show that GRAND($aZ$) is close to
optimal, in the sense of minimizing the number of occupied servers,
when $r$ is large and $a$ is small.

\begin{thm}
\label{th-grand-fluid}
Let $a>0$ be fixed.
Consider a sequence of systems under the GRAND($aZ$) algorithm, indexed by $r$,
and let $\bx^r(\infty)$ denote the random state of the fluid-scaled process
in the stationary regime.
Then, as $r\to\infty$,
$$
\bx^r(\infty) \implies \bx^{*,a}.
$$
\end{thm}

\begin{thm}
\label{th-grand-fluid-convergence}
As $a\downarrow 0$, $\bx^{*,a} \to \cx^*$ and $\bnu^{*,a} \to \ch^*$.
\end{thm}

{\bf Remark 1.} 
Some intuition for the ``entropy-like'' form of the objective function $L^{(a)}$ is as follows.
Speaking informally,
for Theorem~\ref{th-grand-fluid} to hold for some fixed $\bx^{*,a}$, this point $\bx^{*,a}$ must be 
an equilibrium point of the fluid-scaled system under GRAND($aZ$).
(The notion of equilibrium point will be made precise later in the paper.)
As we explain later in Section~\ref{sec-opt-point-interpretation}, 
a fluid-scale equilibrium point under GRAND($aZ$) 
must have a product form; this is shown -- without any consideration of the optimization problem \eqn{eq-opt-grand}-\eqn{eq-cons-laws-grand} --
by interpreting such an equilibrium point 
as the (product-form) stationary distribution of a loss queueing system.
"Entropy-like" objective functions, combined with linear constraints, are known to produce
product-form optimal points; see e.g. papers \cite{zachary2000,kelly1986,kelly1991}, which study the large-scale asymptotics
of a loss system. The fact that entropy maximization leads to a product-form optimal solution is also well-known in the optimization/information theory  literature, see e.g., Example 5.3 in \cite{cvx}. Our choice of "entropy-like" objective $L^{(a)}$, first, results in the product form \eqn{eq-grand-product} of the optimal point $\bx^{*,a}$
and, second, is such that, as $a\to\ 0$, $L^{(a)}$ converges to the linear objective in \eqn{eq-opt}; this motivates its choice.\\
We want to emphasize, however, that while Section~\ref{sec-opt-point-interpretation} provides an interesting interpretation of $\bx^{*,a}$
in terms of a loss system, this interpretation is {\em not} required and {\em not} used in the proofs of our main results.
Also, there is {\em no} close analogy between our main results and those of \cite{zachary2000,kelly1986,kelly1991}: our model, asymptotic regime,
the form and meaning of the optimal point, system dynamics, as well as the main technical idea of the proof of the key Lemma~\ref{lem-fsp-conv} 
(where $L^{(a)}$ is used as a Lyapunov function),
are substantially different
from those in \cite{zachary2000,kelly1986,kelly1991}.

{\bf Remark 2.} Theorem~\ref{th-grand-fluid} is also valid for a closed system with
a fixed number $\rho_i r$ of customers of each type $i$. (Such a closed system,
under different control algorithms, was considered in \cite{St2012,StZh2012}.)
The proof is essentially same as (in fact, simpler than) the one for the open
system in this paper.

\section{Proof of Theorem~\ref{th-grand-fluid-convergence}}
\label{sec-opt-conv}

We prove this theorem first. It is a statement about solutions to the 
optimization problems, and therefore has nothing to do with the system dynamics.


We easily verify that $|L^{(a)}(\bx) - \sum_{\bk} x_{\bk}| \to 0$ as $a\to 0$, uniformly
in $\bx\in \cx$. Therefore, $\bx^{*,a}$ must converge to $\cx^*$.

Consider any sequence $a \to 0$.
We will show that from any subsequence we can choose a further subsequence, along which
we have convergence $\bx^{*,a} \to \bx^*$, $\bnu^{*,a} \to \bbeta^*$,
 where $\bx^* \in \cx^*$ and $\bbeta^* \in \ch^*$ .

Let a subsequence $a\to 0$ be fixed. Since  $\bx^{*,a} \to \cx^*$,
we can and do choose a further subsequence along which $\bx^{*,a} \to \bx^*$ for some fixed $\bx^{*} \in \cx^*$.
Let us show that 
\beql{eq-limsup-1}
\limsup_{a\to 0} \sum_i k_i  \nu^{*,a}_i \le 1, ~~\forall \bk,
\end{equation}
\beql{eq-liminf-0}
\liminf_{a\to 0} \nu^{*,a}_i \ge 0, ~~\forall i.
\end{equation}
If \eqn{eq-limsup-1} would not hold for some $\bk$, then by \eqn{eq-grand-product}
we would have $\limsup x^{*,a}_{\bk} = \infty$ -- a contradiction. Thus, \eqn{eq-limsup-1} holds.
Suppose now that \eqn{eq-liminf-0} does not hold for some $i$, that is 
$\liminf \nu^{*,a}_i < 0$. Pick a $\bk$ with $k_i \ge 1$ and $x^*_{\bk}>0$.
Such a $\bk$ must exist, because $\sum_{\bk} k_i x^*_{\bk} = \rho_i$ (recall that 
$\bx^* \in \cx^*$). Since $x^{*,a}_{\bk} \to x^*_{\bk} \in [0,\rho_i]$,
we see from \eqn{eq-grand-product}
that $\lim \sum_j k_j  \nu^{*,a}_j = 1$. Therefore, 
$$
\limsup \left[\sum_{j\ne i} k_j  \nu^{*,a}_j +
(k_i - 1) \nu^{*,a}_i\right] = 1- \liminf  \nu^{*,a}_i > 1;
$$
 but, this violates \eqn{eq-limsup-1}
for configuration $\bk-\be_i$. Thus, \eqn{eq-liminf-0} holds.

By \eqn{eq-limsup-1}-\eqn{eq-liminf-0}, the sequence of $\bnu^{*,a}$ is bounded.
Then, we choose a further subsequence along which $\bnu^{*,a}$ converges to some $\bbeta^*$.
For the pair $\bx^*$ and $\bbeta^*$, condition \eqn{eq-dual-1} is automatic,
conditions \eqn{eq-dual-111}-\eqn{eq-dual-2} follow from
\eqn{eq-limsup-1}-\eqn{eq-liminf-0}, and condition \eqn{eq-dual-3} follows
from \eqn{eq-grand-product}. Therefore, $\bbeta^* \in \ch^*$. $\Box$

\section{Fluid sample paths under GRAND($aZ$). Proof of Theorem~\ref{th-grand-fluid}}
\label{sec-fsp-dynamics}

In this section, we define fluid sample paths (FSP), 
for the system controlled by GRAND($aZ$)  with $a>0$.
FSPs 
arise as limits of the (fluid-scaled) trajectories $(1/r) \bX^r(\cdot)$
as $r\to\infty$. The definition has a lot in common with the FSP
definitions in \cite{St2012} (and in general is fairly standard),
but it necessarily involves constructions specific to the GRAND($aZ$) discipline
with $a>0$. We \yuan{now} provide details.

Let $\cm$ denote the set of pairs $(\bk,i)$ such that $\bk\in\ck$ and $\bk-\be_i\in\bar\ck$.
Each pair $(\bk,i)$ is associated with the ``edge'' $(\bk-\be_i,\bk)$ connecting configurations
$\bk-\be_i$ and $\bk$; often we refer to this edge as $(\bk,i)$.
By ``arrival along the edge $(\bk,i)$'', we will mean placement of a type $i$ customer
into a server configuration $\bk-\be_i$ to form configuration $\bk$. Similarly, ``departure along the edge $(\bk,i)$'' is a departure of a type-$i$ customer
from a server in configuration $\bk$, which changes its configuration to $\bk-\be_i$.

Without loss of generality, assume that
the Markov process $X^r(\cdot)$ for each $r$ is driven by the
common set of primitive processes, defined as follows.

For each $(\bk,i)\in \cm$, consider an independent
unit-rate Poisson process $\{\Pi_{\bk i}(t), ~t\ge 0\}$, which drives 
departures along edge $(\bk,i)$. Namely, 
let $D^r_{\bk i}(t)$ denote the total 
number of  departures along the edge $(\bk,i)$ in $[0,t]$; then
\beql{eq-driving-dep}
D^r_{\bk i}(t) = \Pi_{\bk i} \left(\int_{\bZero}^t X_{\bk}^r(s) k_i \mu_i ds\right).
\end{equation}
The functional strong law of large numbers (FSLLN) holds:
\beql{eq-flln-poisson-dep}
\frac{1}{r}\Pi_{\bk i}(rt) \to t, ~~~u.o.c., ~~w.p.1.
\end{equation}
For each $i\in \ci$, consider an independent
unit-rate Poisson process $\Pi_{i}(t), ~t\ge 0$, which drives 
exogenous arrivals of type $i$. Namely, 
let $A^r_{i}(t)$ denote the total 
number of type-$i$ arrivals in $[0,t]$, then
\beql{eq-driving-arr}
A^r_{i}(t) = \Pi_{i}(\lambda_i r t).
\end{equation}
Analogously to \eqn{eq-flln-poisson-dep},
\beql{eq-flln-poisson-arr}
\frac{1}{r}\Pi_{i}(rt) \to t, ~~~u.o.c., ~~w.p.1.
\end{equation}
The random placement of new arrivals is constructed as follows.
For each $i\in \ci$, consider an i.i.d. sequence
$\xi_i(1), \xi_i(2), \ldots$ of random variables, uniformly distributed in $[0,1]$.
Denote by $\ck_i \doteq \{\bk\in \bar\ck ~|~ \bk+\be_i \in \bar\ck\}$ the subset of those 
configurations (including zero configuration) which can fit an additional type-$i$ 
customer.
The configurations $\bk\in \ck_i$ are indexed by $1,2,\ldots,|\ck_i|$ 
(in arbitrary fixed order). When the $m$-th (in time) customer 
of type $i$ arrives in the system, it is assigned as follows.
If $X_{(i)}^r=0$, the customer is assigned to an empty server.
Suppose $X_{(i)}^r \ge 1$. Then, the customer
is assigned 
to a server in configuration $\bk'$ indexed by $1$ if 
$$
\xi_i(m) \in [0,X^r_{\bk'}/ X^r_{(i)}],
$$
it is assigned to
a server in configuration $\bk''$ indexed by $2$ if 
$$
\xi_i(m) \in (X^r_{\bk''}/ X^r_{(i)},(X^r_{\bk'}+X^r_{\bk''})/ X^r_{(i)}],
$$
and so on.  Denote
$$
g^r_i(s,\zeta) \doteq \sum_{m=1}^{\lfloor rs \rfloor} I\{\xi_i(m) \le \zeta\},
$$
where $s\ge 0$, $0\le \zeta \le 1$, and $\lfloor \cdot \rfloor$ denotes the integer part of a number. 
Obviously, from the strong law of large numbers
 and the monotonicity of $g^r_i(s,\zeta)$ on both arguments, we 
have the FSLLN
\beql{eq-flln-random}
g^r_i(s,\zeta) \to s\zeta, ~~~\mbox{u.o.c.} ~~~\mbox{w.p.1}
\end{equation}

It is easy and standard to see that, for any $r$, w.p.1,
 the realization of the process $\{\bX^r(t), ~t\ge 0\}$, 
is uniquely determined by the
initial state $\bX^r(0)$
and the realizations of the driving processes
$\Pi_{\bk i}(\cdot)$, $\Pi_{i}(\cdot)$ and $(\xi_i(1), \xi_i(2), \ldots)$.

If we denote by $A^r_{\bk i}(t)$ the total number of arrivals 
allocated along edge $(\bk,i)$ in $[0,t]$, we obviously have
$\sum_{\bk\in \ck_i} A^r_{\bk i}(t)= A^r_{i}(t), ~t\ge 0$, for each $i$

In addition to
$$
x^r_{\bk}(t) = \frac{1}{r} X^r_{\bk}(t),
$$
we introduce other  fluid-scaled 
quantities:
$$
d^r_{\bk i}(t) = \frac{1}{r} D^r_{\bk i}(t),~~~
a^r_{\bk i}(t) = \frac{1}{r} A^r_{\bk i}(t).
$$

A set of locally Lipschitz continuous functions
$[\{x_{\bk}(\cdot),~\bk\in \ck\}, \{d_{\bk i}(\cdot),~(\bk,i)\in \cm\},\{a_{\bk i}(\cdot),~(\bk,i)\in \cm\}]$
on the time interval $[0,\infty)$ we call a {\em fluid sample path} (FSP), if there exist
realizations of the primitive driving processes,
 satisfying conditions \eqn{eq-flln-poisson-dep},\eqn{eq-flln-poisson-arr}
and \eqn{eq-flln-random}
and a fixed subsequence of $r$, along which
\begin{eqnarray}
& [\{x_{\bk}^r(\cdot),~\bk\in \ck\}, \{d_{\bk i}^r(\cdot),~(\bk,i)\in \cm\},\{a_{\bk i}^r(\cdot),~(\bk,i)\in \cm\}]
\to \nonumber \\
& [\{x_{\bk}(\cdot),~\bk\in \ck\}, \{d_{\bk i}(\cdot),~(\bk,i)\in \cm\},\{a_{\bk i}(\cdot),~(\bk,i)\in \cm\}],
~~u.o.c. \label{eq-fsp-def-closed}
\end{eqnarray}

For any FSP, all points $t>0$ are regular (see definition in Section~\ref{subsec-notation}),
except a subset of zero Lebesgue measure. 

\begin{lem}
\label{lem-conv-to-fsp-closed}
Consider a sequence of fluid-scaled processes $\{\bx^r(t),~t\ge 0\}$,
with fixed initial states $\bx^r(0)$ such that $\bx^r(0)\to \bx(0)$.
Then w.p.1, for any subsequence of $r$ there exists a further
subsequence of $r$, along which the convergence \eqn{eq-fsp-def-closed} holds,
with the limit being an FSP.
\end{lem}

{\em Proof} is very standard (see, e.g. \cite{St2012}). 
Essentially, it suffices to observe that, with probability 1, 
on any finite time interval $[0,T]$, $\bx^r(\cdot)$ is uniformly bounded 
for all large $r$. Then, in $[0,T]$, 
the sequences (in $r$)
of functions $d^r_{\bk i}(\cdot)$ and $a^r_{i}(\cdot)$ 
are asymptotically Lipschitz. Namely, for some $C>0$, 
and all $0\le t_1 \le t_2 <\infty$, 
$$
\limsup_r \{d^r_{\bk i}(t_2) - d^r_{\bk i}(t_1)\} \le C(t_2-t_1),
$$
which in turn follows from \eqn{eq-flln-poisson-dep}; and similarly for functions
$a^r_{i}(\cdot)$, and then for $a^r_{\bk i}(\cdot)$ as well.
Using this,  we
easily verify the u.o.c. convergence \eqn{eq-fsp-def-closed}
to some limit with locally Lipschitz components,
along possibly a further subsequence.
We omit details.
$\Box$

For an FSP, at a regular time point $t$, we denote
$v_{\bk i}(t)=(d/dt)a_{\bk i}(t)$ and $w_{\bk i}(t)=(d/dt)d_{\bk i}(t)$.
In other words, $v_{\bk i}(t)$ and $w_{\bk i}(t)$ are the rates of type-$i$ ``fluid''
arrival and departure along edge $(\bk,i)$, respectively. Also
denote: $y_i(t)=\sum_{\bk} k_i x_{\bk}(t)$, $z(t)=\sum_i y_i(t)$, $x_{\bZero}(t)=az(t)$, 
and $x_{(i)}(t) = x_{\bZero} (t) + \sum_{\bk \in \ck: \bk + \be_i \in \ck} x_{\bk}(t)$.

\begin{lem}
\label{lem-fsp-properties-closed-basic}
(i) An FSP satisfies the following properties
at any regular point $t$:
\beql{eq-y-ode}
(d/dt) y_i(t) = \lambda_i - \mu_i y_i(t), ~~\forall i\in \ci,
\end{equation}
\beql{eq-grand-dep-rate}
w_{\bk i}(t)=k_i \mu_i x_{\bk}(t), ~~\forall (\bk,i)\in \cm,
\end{equation}
\beql{eq-grand-arr-rate}
x_{(i)}(t)>0 ~~\mbox{implies}~~
v_{\bk i}(t)=\frac{x_{\bk-\be_i}(t)}{x_{(i)}(t)} \lambda_i, ~~\forall (\bk,i)\in \cm,
\end{equation}
\beql{eq-conserv-rate}
\sum_{\bk:(\bk,i)\in \cm} v_{\bk i}(t) = \lambda_i, ~~\forall i\in \ci,
\end{equation}
\beql{eq-main-diff}
(d/dt) x_{\bk}(t) = \left[\sum_{i:\bk-\be_i\in\bar\ck} v_{\bk i}(t) - \sum_{i:\bk+\be_i\in\bar\ck} v_{\bk+\be_i,i}(t) \right]
                     - \left[\sum_{i:\bk-\be_i\in\bar\ck} w_{\bk i}(t) - \sum_{i:\bk+\be_i\in\bar\ck} w_{\bk+\be_i,i}(t) \right], ~~\forall \bk\in\ck.
\end{equation}
Clearly, \eqn{eq-y-ode} implies
\beql{eq-y-ode2}
y_i(t) = \rho_i + (y_i(0)-\rho_i) e^{-\mu_i t}, ~~t\ge 0, ~~\forall i\in \ci.
\end{equation}
(ii) Moreover, an FSP with $\bx(0)\in \cx$ satisfies the following stronger
conditions:
\beql{eq-y-ode-star}
y_i(t) \equiv \rho_i, ~~\forall i\in \ci,
\end{equation}
\beql{eq-a}
z(t) \equiv 1, ~~ x_{\bZero}(t)\equiv a, ~~ x_{(i)}(t) \ge a, ~\forall i\in \ci;
\end{equation}
at any regular point $t$,
\beql{eq-grand-arr-rate-star}
v_{\bk i}(t)=\frac{x_{\bk-\be_i}(t)}{x_{(i)}(t)} \lambda_i, ~~\forall (\bk,i)\in \cm,
\end{equation}
\beql{eq-conserv-rate-star}
\sum_{\bk:(\bk,i)\in \cm} w_{\bk i}(t) = 
\lambda_i, ~~\forall i\in \ci.
\end{equation}
\end{lem}

{\em Proof.} (i) Given the convergence
\eqn{eq-fsp-def-closed}, \yuan{which defines} an FSP, all the stated properties
except \eqn{eq-grand-arr-rate}, are
nothing  but the limit versions of the flow conservations laws.
Property \eqn{eq-grand-arr-rate} follows from the construction 
of \yuan{the} random assignment, the continuity of $\bx(t)$, and \eqn{eq-flln-random}.
We omit further details. \\
(ii) If $\bx(0)\in \cx$, which implies $y_i(0)=\rho_i$ for each $i$,
property \eqn{eq-y-ode-star} (and then \eqn{eq-a} as well) 
follows from \eqn{eq-y-ode2}. Then, \eqn{eq-grand-arr-rate}
strengthens to \eqn{eq-grand-arr-rate-star}, and 
\eqn{eq-conserv-rate-star} is verified directly using \eqn{eq-grand-dep-rate}.
$\Box$


The key result of this section is the following
\begin{lem}
\label{lem-fsp-conv}
For any FSP with $\bx(0)\in \cx$,
\beql{eq-fsp-conv}
\bx(t) \to \bx^{*,a},
\end{equation}
and the convergence is uniform across all such FSPs.
\end{lem}

{\em Proof.} 
Given that $x_{\bZero}(t) \equiv a$ and $\sum_{\bk} x_{\bk}(t) \le 1$, we have $x_{(i)}(t) \le 1+a$, 
hence $v_{\bk i}(t) \ge x_{\bk}(t) \lambda_i / (1+a)$.
From here, we obtain the following fact: for any $\bk$ and
any $\delta>0$ there exists $\delta_1>0$ such that for all $t\ge \delta$,
 $x_{\bk}(t) \ge \delta_1$. The proof is by 
contradiction. Consider a $\bk$ that is a minimal counterexample; 
necessarily, $\bk \ne \bZero$.
Pick any $\delta>0$ and then the corresponding $\delta_1>0$ such that 
the statement holds for any $\bk'\in \bar\ck$, $\bk'<\bk$. 
(Here $\bk'<\bk$ means that $\bk'\le \bk$ and $\bk'\ne \bk$.)
We observe from \eqn{eq-main-diff} that
 for any regular $t\ge \delta$, $(d/dt) x_{\bk}(t) > \delta_2 >0$ as long as
$x_{\bk}(t) \le \delta_3$, for some positive constants $\delta_2, \delta_3$.
Since this holds for an arbitrarily small $\delta>0$ (with 
$\delta_1, \delta_2, \delta_3$ depending on it), we see that the statement
is true for $\bk$.

In particular, we see that $x_{\bk}(t)>0$ for all $t>0$ and all $\bk$.
Note also that all $t>0$ are regular points (because all $w_{\bk i}$
and $v_{\bk i}$ are bounded continuous in $\bx$).

To prove the lemma, it will suffice to show that:\\
 (a) if $\bx(t) \ne  \bx^{*,a}$ 
and $x_{\bk}(t)>0$ 
for all $\bk \in \ck$, then $(d/dt) L^{(a)}(\bx(t))<0$; and, moreover,\\
(b) the derivative is bounded away from zero as long as $\| \bx(t) -  \bx^{*,a}\|$ 
is bounded away from zero.\\
Let us denote by $\Xi(\bx)$ the derivative $(d/dt) L^{(a)}(\bx(t))$ at a given point
$\bx(t)=\bx$; in the rest of the proof we study the function $\Xi(\bx)$ on $\cx$, 
and therefore drop the time index $t$. Suppose all components $x_{\bk}>0$.
From \eqn{eq-grand-dep-rate}, \eqn{eq-conserv-rate},
\eqn{eq-grand-arr-rate-star}, \yuan{and}
\eqn{eq-conserv-rate-star}, we have:
$$
w_{\bk i} = k_i \mu_i x_{\bk} \sum_{\bk': (\bk',i) \in \cm} \frac{x_{\bk'-\be_i}}{x_{(i)}},
$$
$$
v_{\bk' i} = \frac{x_{\bk'-\be_i}}{x_{(i)}} \sum_{\bk: (\bk,i) \in \cm} k_i \mu_i x_{\bk}.
$$
We can use the following interpretation of these expressions for $w_{\bk i}$ and $v_{\bk' i}$: 
for any ordered pair of edges $(\bk,i)$ and
$(\bk',i)$, we can assume that the part $k_i \mu_i x_{\bk} x_{\bk'-\be_i}/x_{(i)}$
of the total departure rate $k_i \mu_i x_{\bk}$ along $(\bk,i)$ is ``allocated back''
as a part of the arrival rate along $(\bk',i)$. 
 (This interpretation is not needed, in principle. But, we think it helps to see
  the intuition behind the formal expressions that follow.)
The contribution of these ``coupled'' departure/arrival
rates for the ordered pair of edges $(\bk,i)$ and
$(\bk',i)$ into the derivative $\Xi(\bx)$ is 
(after simple manipulation):
$$
\xi_{\bk,\bk',i} = (1/b)\left[\log (k'_i x_{\bk-\be_i} x_{\bk'}) - \log (k_i x_{\bk} x_{\bk'-\be_i}) \right]
\frac{k_i \mu_i x_{\bk} x_{\bk'-\be_i}}{x_{(i)}}.
$$
We remark that this expression is valid even when either $\bk-\be_i=\bZero$ or $\bk'-\be_i=\bZero$.
This is because 
$x_{\bZero}(t) = a$ when $\bx\in \cx$, and therefore by convention \eqn{eq-formal-deriv-zero},
formula \eqn{eq-L-partial} is valid for all $\bk\in \bar\ck$. We have:
$$
\xi_{\bk,\bk',i} + \xi_{\bk',\bk,i} 
= (1/b)(\mu_i/x_{(i)}) [\log (k'_i x_{\bk-\be_i} x_{\bk'}) - \log (k_i x_{\bk} x_{\bk'-\be_i}) ] 
[k_i x_{\bk} x_{\bk'-\be_i} - k'_i x_{\bk-\be_i} x_{\bk'}] \le 0,
$$
and the inequality is strict unless $k'_i x_{\bk-\be_i} x_{\bk'}=k_i x_{\bk} x_{\bk'-\be_i}$.
We obtain
\beql{eq-L-deriv}
\Xi(\bx) = \sum_i \sum_{\bk,\bk'} [\xi_{\bk,\bk',i} + \xi_{\bk',\bk,i}].
\end{equation}
Therefore, 
$\Xi(\bx)<0$ unless $\bx$ has a product form representation \eqn{eq-grand-product},
which in turn is equivalent to $\bx=\bx^{*,a}$.

So far the function $\Xi(\bx)$ in \eqn{eq-L-deriv} was defined for $\bx\in\cx$
with all $x_{\bk}>0$. Let us adopt a convention that $\Xi(\bx)=-\infty$
for $\bx\in\cx$ with at least one $x_{\bk}=0$. Then,
it is easy to verify that $\Xi(\bx)$ is continuous on the entire set $\cx$.

It remains to show that
for any $\delta_2 > 0$ there exists $\delta_3 > 0$ such that conditions
$\bx\in \cx$ and
$L^{(a)}(\bx) - L^{(a)}(\bx^{*,a}) \ge \delta_2$ imply
$\Xi(\bx) \le -\delta_3$. This is indeed true, because otherwise
there would exist $\bx\in \cx$, $\bx\ne \bx^{*,a}$, such that 
 $\Xi(\bx) = 0$, which is, again, equivalent to $\bx = \bx^{*,a}$.
$\Box$

From Lemma~\ref{lem-fsp-conv} we easily obtain Theorem~\ref{th-grand-fluid}.
The argument is essentially same as in the conclusion 
of the proof of Theorem 2 in \cite{St2012}.

{\em Proof of Theorem~\ref{th-grand-fluid}.} 
We fix $\epsilon>0$ and choose $T>0$ large enough so that, for any FSP 
with $\bx(0)\in \cx$
we have $\|\bx(T)-\bx^{*,a}\| \le \epsilon$.
We claim that for any $\delta_1>0$ there exists a $\delta>0$ such that,
uniformly for all sufficiently large $r$ and all initial 
states $\bx^r(0)$ such that $d(\bx^r(0),\cx)\le \delta$, we have
\beql{eq-claim1}
\pr\{\|\bx^r(T)-\bx^{*,a}\|<2\epsilon\} > 1-\delta_1.
\end{equation}
This claim is true, because 
otherwise we would have a sequence of fixed initial states
$\bx^r(0)\to \bx(0) \in \cx$, for which
$$ 
\pr\{\|\bx^r(T)-\bx^{*,a}\|<2\epsilon\} \le 1-\delta_1;
$$ 
this in turn is impossible
because w.p.1 we can always choose a subsequence of $r$ along which the 
u.o.c convergence to an FSP with initial state $\bx(0)$ holds.
Claim \eqn{eq-claim1} implies the result,
because $\epsilon$ can be arbitrarily small
and (according to \eqn{eq-cons-laws}) the stationary version of the
process is such that $d(\bx^r(0),\cx)\le \delta$ with 
\yuan{a probability that is arbitrarily close to 1}, for all large $r$.
$\Box$

Although not required for the proof of Theorem~\ref{th-grand-fluid},
we note that the following generalization of Lemma~\ref{lem-fsp-conv} holds,
which is of independent interest, because it describes \yuan{the} FSP dynamics for
arbitrary initial states, not necessarily $\bx(0)\in \cx$. 

\begin{lem}
\label{lem-fsp-conv-gen}
For any compact $A \in \R_+^{|\ck|}$, the convergence
\beql{eq-fsp-conv-gen}
\bx(t) \to \bx^{*,a}
\end{equation}
holds uniformly in all FSPs with $\bx(0)\in A$.
\end{lem}

{\em Proof} is a slight generalization of that of Lemma~\ref{lem-fsp-conv}.
From \eqn{eq-y-ode2} we conclude that, starting any fixed time $\tau>0$,
$0< a_1 \le x_{\bZero}(t) \le x_{(i)}(t) \le a_2 < \infty$, for some constants $a_1,a_2$, uniformly 
on all FSPs with $\bx(0)\in A$. (This ``replaces'' condition
$a = x_{\bZero}(t) \le x_{(i)}(t) \le 1+a$.) Then, if we reset the initial time
from $0$ to $\tau$, we obtain (the same way) this fact: 
for any $\bk$ and
any $\delta>0$ there exists $\delta_1>0$ such that for all $t\ge \tau+\delta$,
 $x_{\bk}(t) \ge \delta_1$. In particular, all points $t>0$ are regular,
with all $x_{\bk}(t)>0$. 

At any (regular) point $t>0$, we have the following general expression:
\beql{eq-L-deriv-gen}
(d/dt) L^{(a)}(\bx(t)) = \Xi(\bx(t)),
\end{equation}
where
\beql{eq-L-deriv-gen-Xi}
\Xi(\bx) \doteq \sum_{(\bk,i)\in\cm} [v_{\bk i}-w_{\bk i}] [f(\bk) - f(\bk-\be_i)I\{\bk-\be_i\ne 0\}], \mbox{and}
\end{equation}
$$
f(\bk) = (\partial / \partial x_{\bk}) L^{(a)}(x) =(1/b)  \log [x_{\bk} c_{\bk} / a ].
$$
(Here $I\{\cdot\}$ is the indicator function.)
Expression \eqn{eq-L-deriv-gen-Xi} is more general than \eqn{eq-L-deriv},
which was for the special case $\bx\in \cx$. It is well-defined when all $x_{\bk}>0$.
We extend the definition to those $\bx$ with some $x_{\bk}=0$, in which case
we set $\Xi(\bx)=-\infty$.
Fix any $\epsilon>0$ such that $\epsilon < \min_i \rho_i$.
Consider set $B\doteq \{\bx\in R_+^{|\ck|} ~|~ |y_i-\rho_i| \le \epsilon, ~\forall i\}$.
(Clearly, $\cx \subset B$;
$\bx(t)\in B$ after a finite time, uniformly on all FSPs with $\bx(0)\in A$.)
Using \eqn{eq-L-deriv-gen-Xi},
it is not difficult to verify that $\Xi(\bx)$ is continuous on $B$.

Observe that since $\bx(0)$ is bounded, then $\bx(t)$ (and therefore $L(\bx(t))$)
is uniformly bounded.

It remains to prove the following assertion:
for any $\delta_2 > 0$ there exist $\delta_3 > 0$ 
and a sufficiently large $T>0$
such that condition
$L^{(a)}(\bx(t)) - L^{(a)}(\bx^{*,a}) \ge \delta_2$ implies $(d/dt)L^{(a)}(\bx(t)) \le -\delta_3$. 
Indeed, for a given $\delta_2 > 0$ we first choose $\delta_3 > 0$ 
the same way as in the analogous
assertion at the end of the proof of Lemma~\ref{lem-fsp-conv}.
Fix a small $\delta_4>0$ and then choose $T>0$ large enough so that (by \eqn{eq-y-ode2})
$d(\bx(t),\cx) < \delta_4$. Using the continuity of $\Xi$ on $B$,
it is easy to check that if $\delta_4$ is small enough,
the assertion must hold with $\delta_3$ rechosen to be, say, $\delta_3/2$.
We omit further details.
$\Box$

\section{Probabilistic interpretation of the optimal point $\bx^{*,a}$.}
\label{sec-opt-point-interpretation}

Lemma~\ref{lem-fsp-conv} implies, in particular, that
$\bx(t)\equiv \bx^{*,a}$ is a (unique) stationary FSP. 
That is, $\bx^{*,a}$ is an equilibrium point of the fluid system (i.e., the dynamic 
system described by FSPs).
We now show that {\em if a point  $\bx^{**}$
is equilibrium for the fluid system, this fact alone implies that 
$\bx^{**}$ must have a product form, and therefore $\bx^{**}=\bx^{*,a}$.}
More precisely, the argument below shows that an equilibrium point can be interpreted 
as a (rescaled) unique stationary distribution of a corresponding loss system.

If $\bx(t)\equiv \bx^{**}$ is an FSP, then $\bx^{**}$ satisfies the stationary version
of \eqn{eq-main-diff}, namely,
$$
\sum_{i:\bk-\be_i\in\bar\ck} \hat\lambda_i x_{\bk-\be_i}^{**}
 - \sum_{i:\bk+\be_i\in\bar\ck} \hat\lambda_i x_{\bk}^{**}  -
$$
\beql{eq-main-diff-stat}
\sum_{i:\bk-\be_i\in\bar\ck} k_i \mu_i x_{\bk}^{**} + \sum_{i:\bk+\be_i\in\bar\ck} (k_i+1) \mu_i x_{\bk+\be_i}^{**} = 0, 
~~\forall \bk\in\ck,
\end{equation}
where $x_{\bZero}^{**} = a$ and we denoted $\hat\lambda_i = \lambda_i/x_{(i)}^{**}$.
We observe that equations \eqn{eq-main-diff-stat} are exactly the same as
balance equations for the steady-state probabilities of the following loss
queueing system.  The system has exogenous Poisson input customer flows with
 rates $\hat\lambda_i$. The service time of a type-$i$ customer is an independent and
exponentially distributed random variable with mean $1/\mu_i$.
The feasible system states are given by configurations $\bk\in\bar\ck$.
If an arriving type-$i$ customer finds the system in a state $\bk$ such that
$\bk+\be_i\not\in \bar\ck$, this customer is lost. It is well known 
that
the stationary distribution of this
system is unique and has product form
$$
x_{\bk}^{**} = x_{\bZero}^{**} \prod_i \frac{1}{k_i !} \left(\frac{\hat\lambda_i}{\mu_i}\right)^{k_i}
= x_{\bZero}^{**} \frac{1}{c_{\bk}} \prod_i  \left(\frac{\hat\lambda_i}{\mu_i}\right)^{k_i}.
$$
For $\bx^{**}$, instead of the normalization condition
$\sum_{\bk\in \bar\ck} x_{\bk}^{**} = 1$ we have a ``boundary condition''
$x_{\bZero}^{**} = a$, so finally we can write
\beql{eq-grand-product2}
x_{\bk}^{**} 
= a \frac{1}{c_{\bk}} \prod_i  \left(\frac{\hat\lambda_i}{\mu_i}\right)^{k_i}.
\end{equation}
Since \eqn{eq-grand-product2} can be rewritten in form \eqn{eq-grand-product}, this implies that $\bx^{**}=\bx^{*,a}$.

Comparing \eqn{eq-grand-product2} and \eqn{eq-grand-product}, we obtain 
the relation between $\hat \lambda_i$ (and then $x_{(i)}^{*,a}=\lambda_i/\hat\lambda_i$) and $\nu_i^{*,a}$:
\beql{eq-hat-lambda}
\hat \lambda_i = \mu_i a^{-\nu_i^{*,a}}.
\end{equation}
We note that the loss system described above is, of course, not just 
an abstract construction. If we consider the system under GRAND($aZ$), with large $r$, in 
steady state, the random evolution of any server that is either occupied or one of the zero servers is close to that
of the described loss system. Namely, this server ``experiences'' a type-$i$ arrival
process that is close to a Poisson process of constant rate
$ \hat\lambda_i = \lambda_i/x_{(i)}^{*,a}=\mu_i a^{-\nu_i^{*,a}}$, as long as the server has room
for at least one more type $i$ customer (and as long as the server remains occupied, or belongs 
to one of the zero servers).

{\bf Remark 3.} It is further well known that the stationary distribution of states $\bk$ of a loss queueing system,
described in this section, is insensitive to detailed service time distributions, as long as mean service times $1/\mu_i$ are fixed.
(See e.g. \cite{insensitivityZ,insensitivityB}.) Given this fact, the argument in this section can be repeated (essentially as is)
to show the following. Suppose we consider the model and the asymptotic regime of this paper, generalized to allow
arbitrary service time distributions, with mean values $1/\mu_i$. Consider FSPs for such a system, which are defined more generally,
because a system state will now include customer residual (or elapsed) service times. 
Then, any equilibrium point of the fluid system
 is equal (up to scaling by a positive constant) 
to a stationary distribution of the corresponding loss system with input rates \eqn{eq-hat-lambda}.
The latter stationary distribution is unique.
So, we obtain that {\em the equilibrium point of the fluid system is unique
and is equal (up to scaling by a positive constant) 
to the unique stationary distribution of the corresponding loss system with input rates \eqn{eq-hat-lambda}.}
In particular, the 
$\bx$-projection of the equilibrium point is still equal to $\bx^{*,a}$, as defined in this paper.
This strongly suggests that our Theorem \ref{th-grand-fluid} holds as is for the system with general service time distributions.

\section{Simulations}
\label{sec-sim}

In previous sections, we proved the asymptotic optimality of GRAND($aZ$), as 
the system scale (customer arrival rates) grows to infinity and parameter
$a \rightarrow 0$. 
This suggests that  GRAND algorithms with fewer zero-servers 
result in smaller numbers of occupied servers in steady state, hence better performance. 
On the other hand, intuitively, more zero-servers create greater ``safety stocks''
of servers in many different configurations, which should help the system move
to its steady state faster. In other words, although a larger number of zero-servers
results in greater steady-state inefficiency (because, roughly speaking, $\bx^{*,a}$ with larger $a$ is farer from $\cx^*$),
one might expect a faster convergence, and perhaps a better performance at least 
in some time intervals during a transition to steady state.

Motivated by these considerations, in this section we conduct simulation experiments 
to compare different versions of a GRAND algorithm under different settings.  
We will identify a specific instance of GRAND by appending in parentheses the function giving the number of zero-servers.
(Definition \ref{df:grand-az} of GRAND($aZ$) is consistent with this convention.)
In particular, the following instances are considered:
\begin{itemize}
\item[(i)] GRAND($aZ$) for $a>0$; 
\item[(ii)] GRAND($c$), for some fixed non-negative integer constant $c$
(i.e., $X_{\bZero}(t) \equiv c$ for all time $t$).
\end{itemize}

\begin{figure}
\centering
\subfigure[Vector packing; feasible configurations are $\{(k_1, k_2) : 2k_1 + 3k_2 \leq 15\}$.]{\includegraphics[width=7cm]{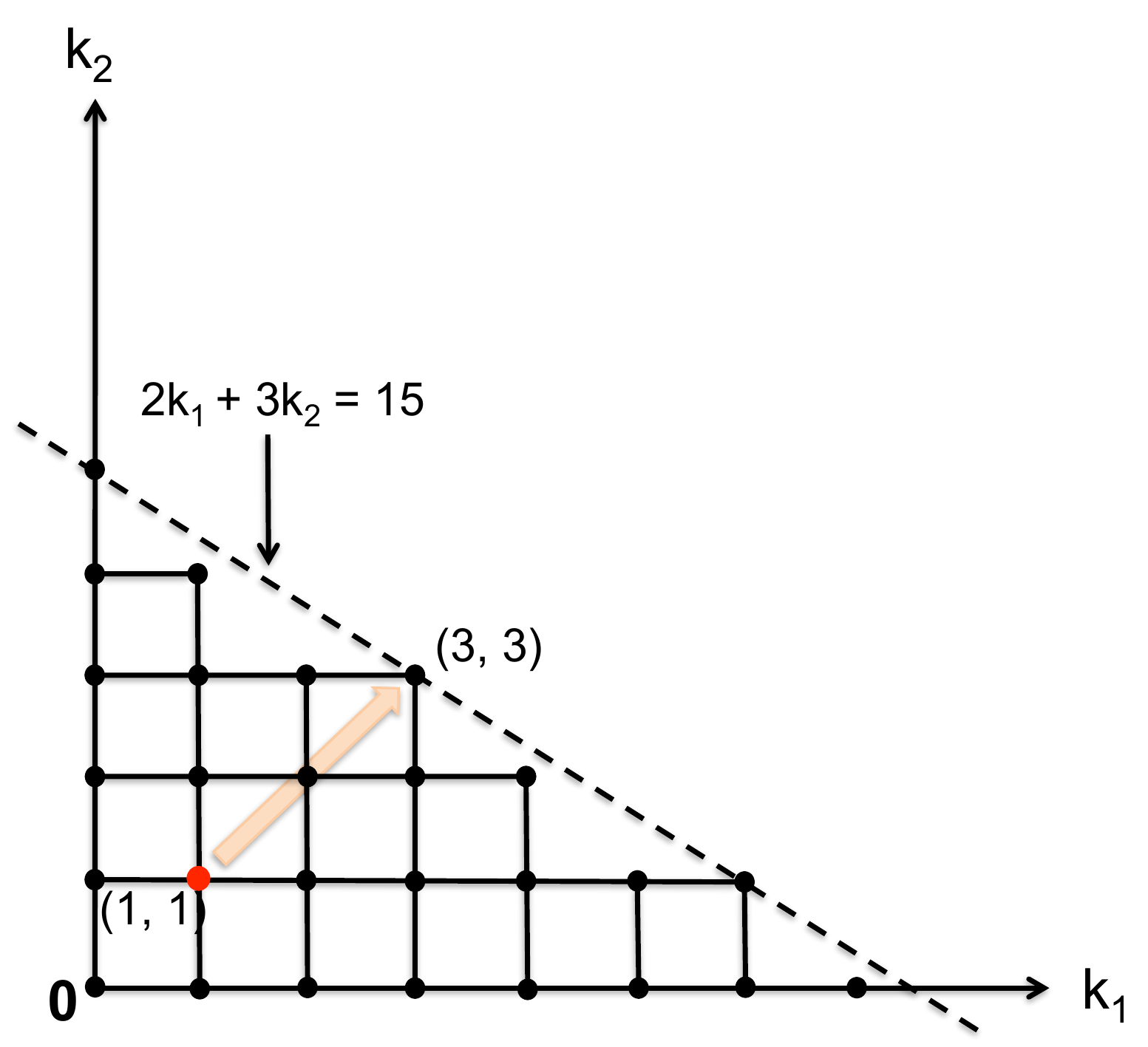}\label{fig:vp}}
\hfill
\subfigure[Not vector packing; maximal configurations are $\{(8, 1), (3, 3), (1, 8)\}$.]{\includegraphics[width=7cm]{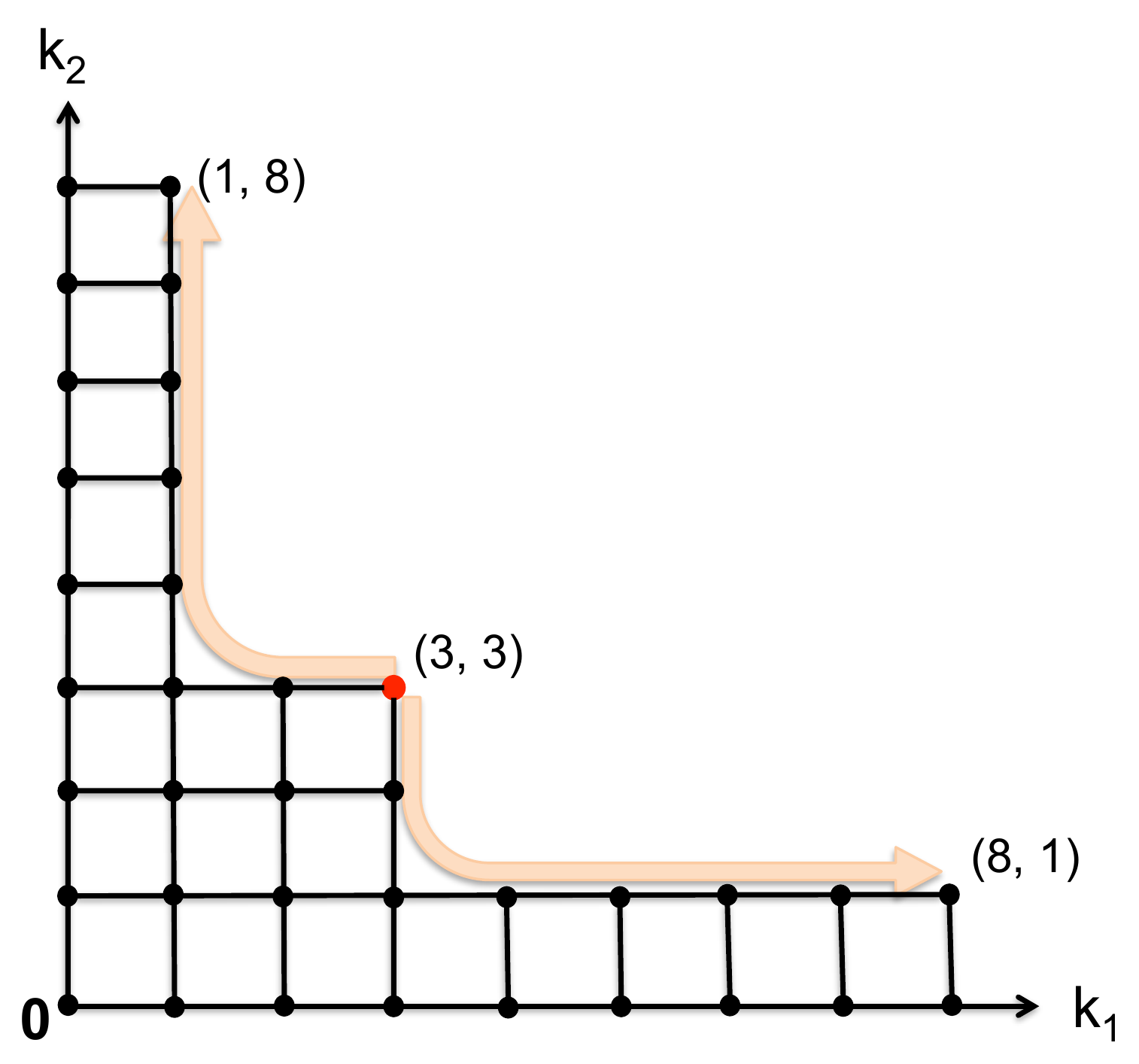}\label{fig:non-vp}}
\caption{Two systems with different sets of feasible configurations.}
\end{figure}

\paragraph{Simulation setup.}  
We consider the performance and convergence times of GRAND algorithms 
in two different systems, illustrated in Figures \ref{fig:vp} and \ref{fig:non-vp}. 
From now on we call these systems  
\ref{fig:vp} and \ref{fig:non-vp} respectively. In \ref{fig:vp} 
the packing constraints are of vector-packing type (cf. \cite{St2012} for the exact definition).
In this system, each server is a one-dimensional bin of size $15$.
There are $2$ types of arriving customers. A type-$1$ customer 
occupies $2$ units of space when in service, and 
a type-$2$ customer, $3$ units of space. As such, 
the set $\bar{\ck}$ of feasible configurations is given by 
$\{(k_1, k_2) \in \mathbb{Z}_+^2 : 2k_1+3k_2 \leq 15\}$ (cf. Figure \ref{fig:vp}). 
In \ref{fig:non-vp} there are also $2$ types 
of arriving customers, but packing constraints are {\em not} of vector-packing type.
  The set $\bar{\ck}$ of feasible server configurations 
is characterized by the maximal server configurations 
$(8, 1), (3, 3)$ and $(1, 8)$ (cf. Figure \ref{fig:non-vp}). 

In \emph{both} systems, the system parameters are given as follows. 
The service time of any customer is independently 
and exponentially distributed with mean $1$, 
so that $\mu_1 = \mu_2 = 1$. The normalized arrival rates 
are $\lambda_1 = \lambda_2 = 1/2$, so that in a system of scale $r$, 
the actual arrival rates are $\Lambda_1 = \Lambda_2 = r/2$. 
For a system of scale $r$, in \ref{fig:vp}, 
the initial system state is always set to be $X_{(1, 1)}(0) = \lfloor r/2 \rfloor$ 
and $X_{\bk}(0) = 0$ for all $\bk \neq (1, 1)$; and in \ref{fig:non-vp}, 
$X_{(3, 3)}(0) = \lfloor r/6 \rfloor$ and $X_{\bk} = 0$ for all $\bk \neq (3, 3)$. 
Note that in both systems the initial total number of customers of each type is equal to its 
average number in steady state; that is, the fluid-scaled initial state is within $\cx$.

It is clear that in \ref{fig:vp}, the fluid-scale optimal solution is 
$x_{(3,3)} = 1/6$ and $x_{\bk} = 0$ for all $\bk \neq (3, 3)$, so that 
$\sum_{\bk \in\ck} x_{\bk}^*$, the optimal total (fluid-scale) number of occupied servers, 
is $1/6$.
In \ref{fig:non-vp}, the fluid-scale optimal solution 
is $x_{(8, 1)} = 1/18$, $x_{(1, 8)} = 1/18$ and $x_{\bk} = 0$ for all $\bk \neq (1, 8), (8, 1)$; 
hence $\sum_{\bk \in \ck} x_{\bk}^*$, the optimal total (fluid-scale) number 
of occupied servers is $1/9$. 
Here we remark that system parameters are chosen 
to simplify structures of the optimal solutions. 
However, observations from simulations 
in systems considered here are representative of those in more general systems. 
We now describe the simulation results in detail. 

\begin{center}
\begin{figure}
\includegraphics[scale=0.4]{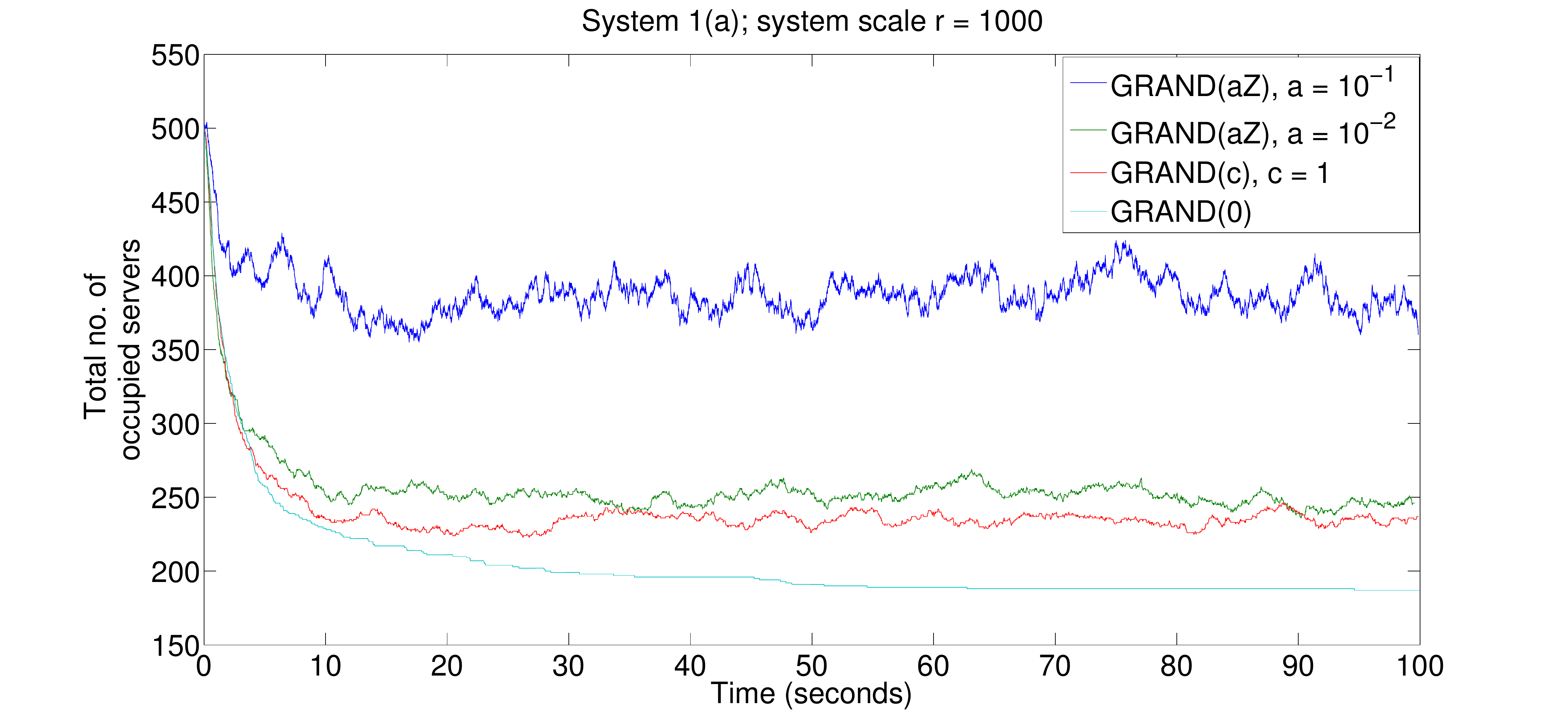}
\caption{System \ref{fig:vp}. The system scale is $r = 1000$, and 
the initial condition is $X_{(1,1)}(0) = r/2 = 500$, and $X_{\bk}(0) = 0$ for $\bk \neq (1, 1)$.}
\label{fig:compare-vp-1000}
\end{figure}
\end{center}

\begin{center}
\begin{figure}
\includegraphics[scale=0.4]{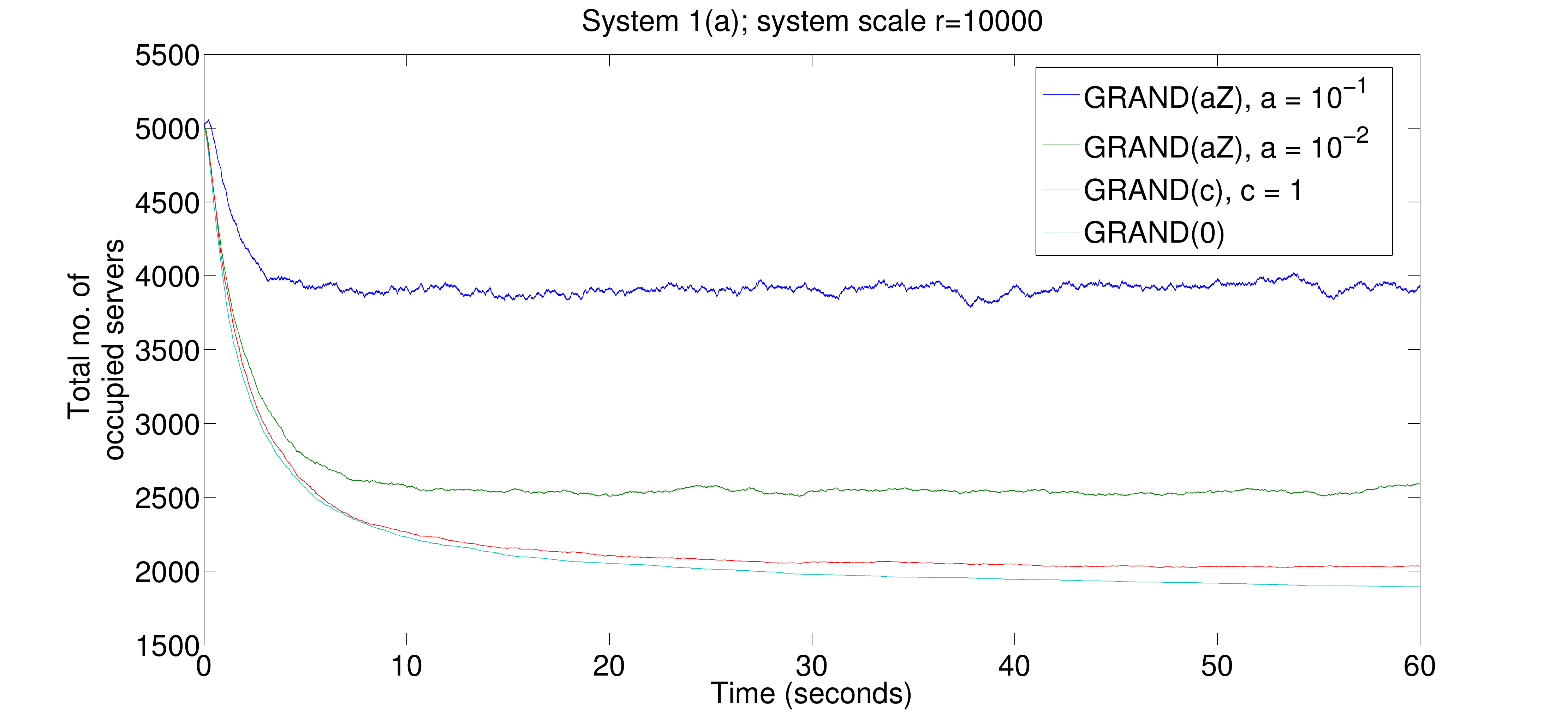}
\caption{System \ref{fig:vp}. The system scale is $r = 10000$, and 
the initial condition is $X_{(1,1)}(0) = r/2 = 5000$, $X_{\bk}(0) = 0$ for $\bk \neq (1, 1)$.}
\label{fig:compare-vp-10000}
\end{figure}
\end{center}

\begin{center}
\begin{figure}
\includegraphics[scale = 0.4]{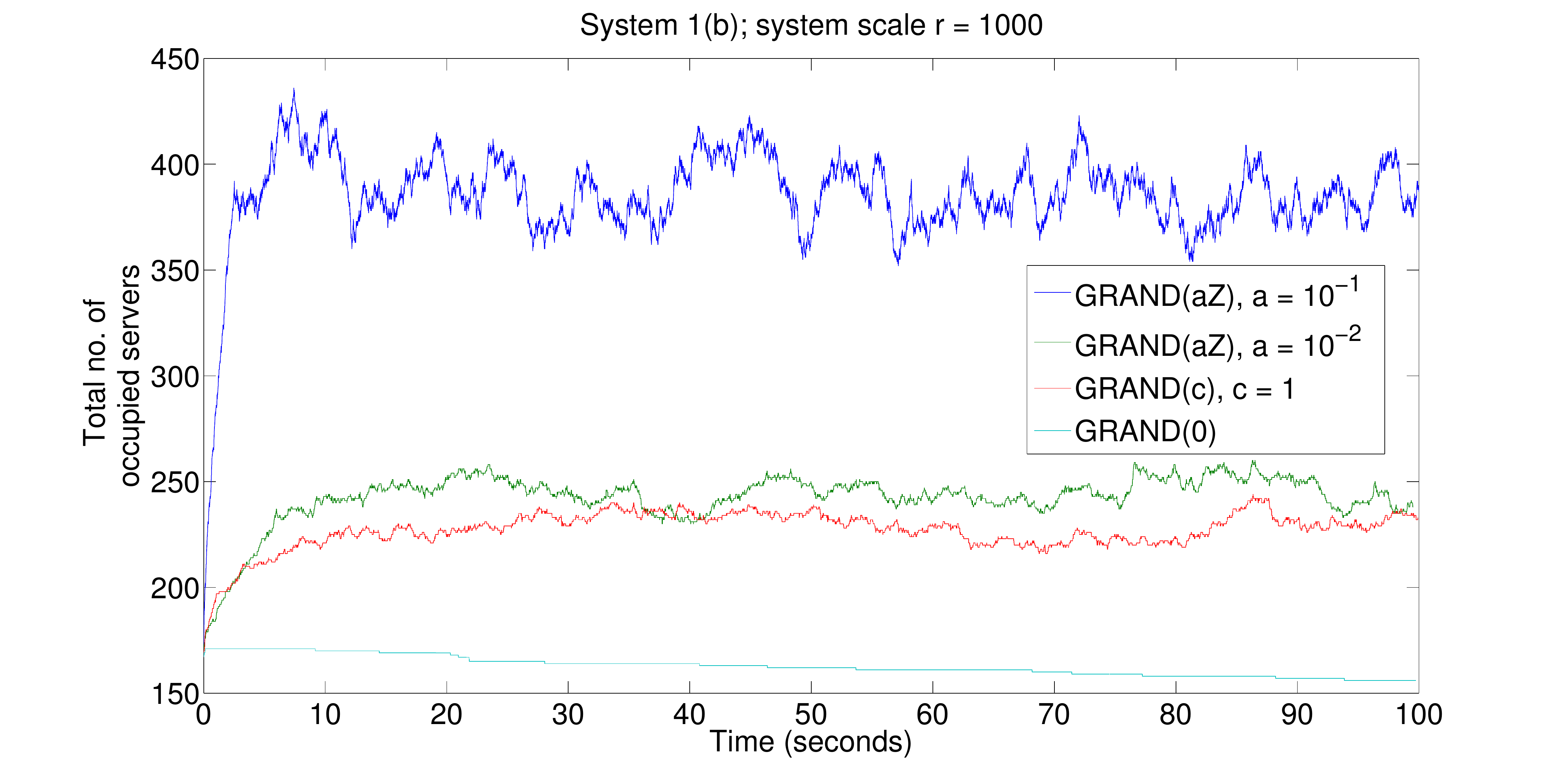}
\caption{System \ref{fig:non-vp}. The system scale is $r = 1000$, 
and the initial condition is $X_{(3,3)}(0) = \lfloor r/6 \rfloor = 166$ 
and $X_{\bk}(0) = 0$ for $\bk \neq (3, 3)$.}
\label{fig:compare-corner-1000}
\end{figure}
\end{center}

\begin{center}
\begin{figure}
\includegraphics[scale=0.4]{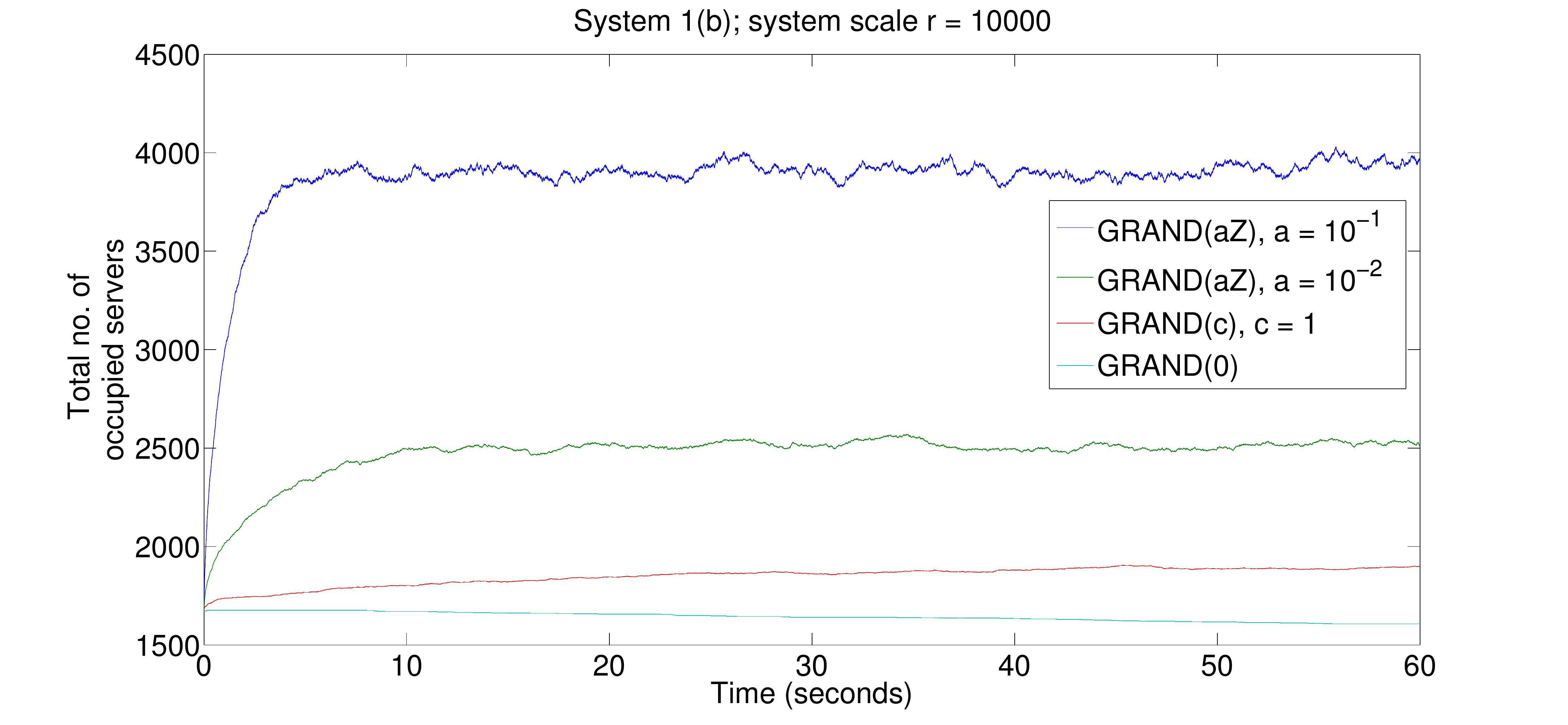}
\caption{System \ref{fig:non-vp}. The system scale is $r = 10000$, and 
the initial condition is $X_{(3,3)(0)} = \lfloor r/6 \rfloor = 1666$ and $X_{\bk} = 0$ for $\bk \neq (3, 3)$.}
\label{fig:compare-corner-10000}
\end{figure}
\end{center}

\paragraph{Observations from discrete-event simulations.} 
In Figures \ref{fig:compare-vp-1000} 
through \ref{fig:compare-corner-10000}, 
we simulate the evolutions 
of systems \ref{fig:vp} and \ref{fig:non-vp} at scales $r = 1000$ and $r = 10000$. 
Table \ref{tab:summary} is a simple summary of the settings. 

The simulated algorithms are GRAND($aZ$) with $a = 0.1$ and $a = 0.01$, 
GRAND($c$) with $c=1$, and GRAND($0$). We have plotted the sample-path trajectories of 
the total number of occupied servers. These plots confirm the tradeoff between 
convergence time and performance: 
as we decrease the number of zero-servers, 
the algorithm performance improves, 
but the convergence time increases. 
For example, in Figure \ref{fig:compare-vp-10000}, 
GRAND($aZ$) with $a = 0.1$ converges to about 4000 occupied servers at around $t = 4$, 
and GRAND($aZ$) with $a = 0.01$ converges to about 2500 occupied servers at around $t = 10$.

However, in these plots, an increased number of zero-servers does not seem to help with performance 
at all, as measured by the total number of occupied servers.
Namely, GRAND($0$) produces the least number 
of occupied servers at almost all times, 
making its slower convergence to steady state seem irrelevant. 
In fact, at scales $r =1000$ and $r = 10000$ 
and in some other systems where we simulate the GRAND algorithms, 
we do not see any performance benefit of using a positive number of zero-servers either.

\paragraph{Observations from fluid-limit dynamics.} 
The above observations motivate us to consider the 
system dynamics under GRAND($aZ$) algorithms, with different $a > 0$, 
in the fluid limit. The fluid-limit trajectories are described by FSPs.
(See Lemma ~\ref{lem-conv-to-fsp-closed}.) We restrict the initial conditions $\bx(0)$ 
to the set $\cx$, i.e., $\bx(0) \in \cx$, so that
FSPs satisfy equations/conditions \eqref{eq-y-ode} -- \eqref{eq-conserv-rate-star}
in Lemma \ref{lem-fsp-properties-closed-basic}. 
In our case, it is quite easy to see that given an initial condition $\bx(0) \in \cx$, 
a trajectory satisfying \eqref{eq-y-ode} -- \eqref{eq-conserv-rate-star} is unique,
and therefore it is the unique FSP, describing precisely 
the system dynamics in the fluid limit. 
This is because at $t = 0$, right-derivatives are well-defined, and for all times 
$t > 0$, $x_{\bk} (t) > 0$ for all $\bk$ (see the proof of Lemma \ref{lem-fsp-conv}) and derivatives 
are continuous and well-defined. Therefore, we can compute the FSP as the unique solution to an ODE.
Our observations are as follows.

\begin{center}
\begin{figure}
\includegraphics[scale=0.4]{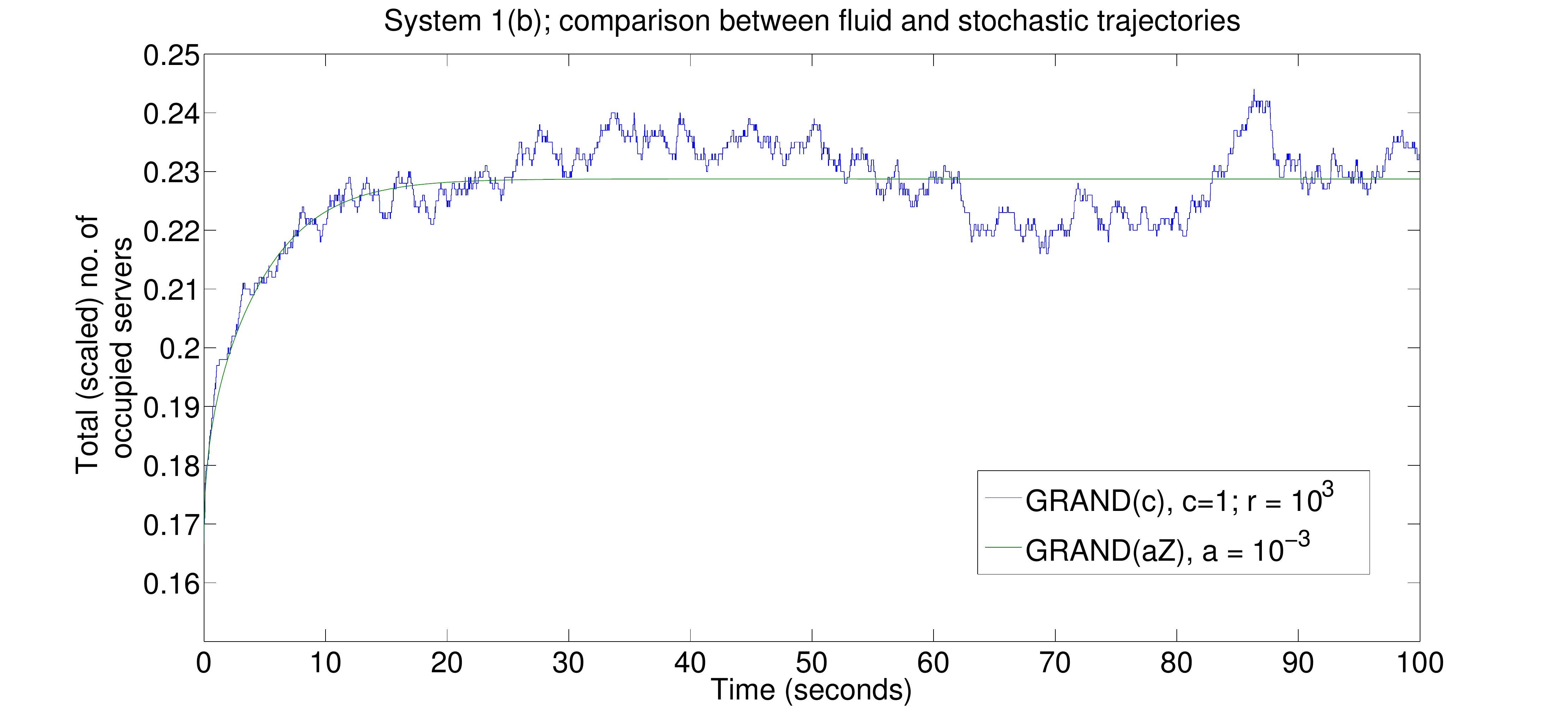}
\caption{System \ref{fig:non-vp}. Comparison between trajectories 
in the fluid limit and of systems at finite scale.}
\label{fig:fluid-stoc}
\end{figure}
\end{center}

First, the fluid-limit dynamics is a reasonable approximation to 
the actual system dynamics, even at moderate scale. 
Figure \ref{fig:fluid-stoc} plots the evolution of the system \ref{fig:non-vp} 
at scale $r = 1000$ under GRAND($c$), $c=1$, 
and compare it with that in the fluid limit under GRAND($aZ$) with $a = 10^{-3}$. 
For the system with scale $r$, we have normalized the 
total number of occupied servers by $r$. 
The value $a=10^{-3}$ is chosen so that, 
 in the system of scale $r = 1000$,
$aZ \approx ar = 1$. 
We see that the trajectories agree reasonably well. 

\begin{center}
\begin{figure}
\includegraphics[scale=0.4]{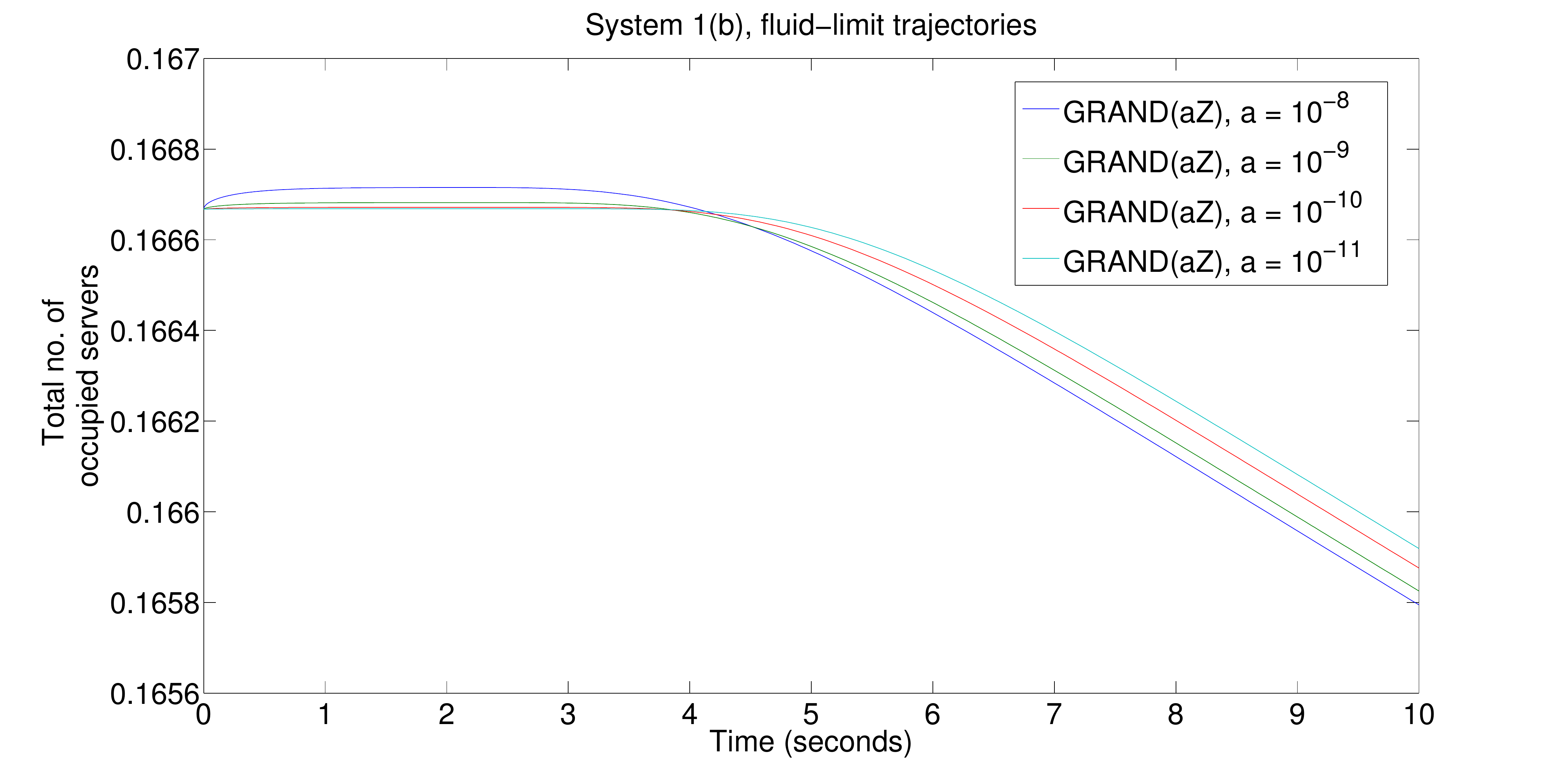}
\caption{Comparison of GRAND($aZ$) for different $a$, in system \ref{fig:non-vp}. 
$x_{(3, 3)}(0) = 1/6$, and $x_{\bk}(0) = 0$ for $\bk \neq (3, 3)$.}
\label{fig:compare-fluid}
\end{figure}
\end{center}

\begin{center}
\begin{figure}
\includegraphics[scale=0.4]{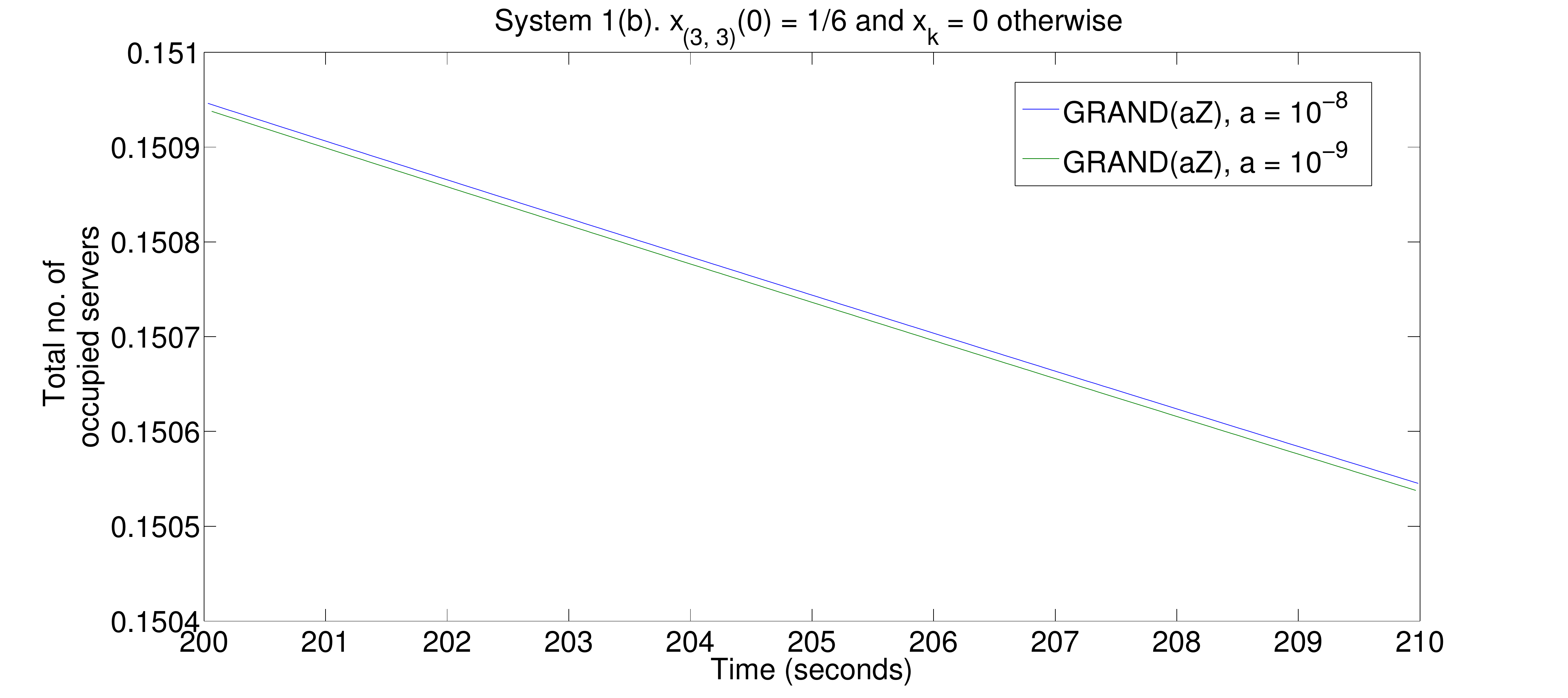}
\caption{Comparison of GRAND($aZ$) between $a = 10^{-8}$ and $a = 10^{-9}$, 
in system \ref{fig:non-vp}. 
$x_{(3, 3)}(0) = 1/6$, and $x_{\bk}(0) = 0$ for $\bk \neq (3, 3)$.}
\label{fig:compare-fluid-large-t}
\end{figure}
\end{center}

Second, in Figure \ref{fig:compare-fluid} we do see 
some performance benefit of a larger number of zero-servers. 
For system \ref{fig:non-vp}, Figure \ref{fig:compare-fluid} plots 
the number of occupied servers  for FSPs
 under GRAND($aZ$) with different $a$.
Compare, for example, GRAND($aZ$) with $a = 10^{-8}$ 
and GRAND($aZ$) with $a = 10^{-9}$. Between time $0$ 
and $4$, GRAND($aZ$) with $a = 10^{-9}$ produces a smaller 
number of occupied servers; between time $5$ and $10$, GRAND($aZ$) 
with $a = 10^{-8}$ produces a smaller number. 
Since $10^{-9} < 10^{-8}$, we also expect 
GRAND($aZ$) with $a = 10^{-9}$ to perform better 
than GRAND($aZ$) with $a = 10^{-8}$ at sufficiently large times. 
This is captured in Figure \ref{fig:compare-fluid-large-t}; 
we did not plot the whole trajectories 
since these performance differences are quite small 
and would not be very visible in a plot at longer time scales. 
We observe that for system \ref{fig:non-vp}, 
we only see the performance benefit of an increased number of zero-servers 
for the values of $a$ that are very small, e.g., $a \leq 10^{-8}$. 
This suggests that the performance benefit (if any)
of an increased number of zero-servers in an actual system can only be observed 
at a very large scale, e.g., when $r \geq 10^8$. 
Thus, in principle, the performance benefit 
of a larger number of zero-servers does exist, 
but it might not be seen in the systems of practical interest such as data centers, 
if they are of moderate scale.  

\begin{center}
\begin{figure}
\includegraphics[scale=0.4]{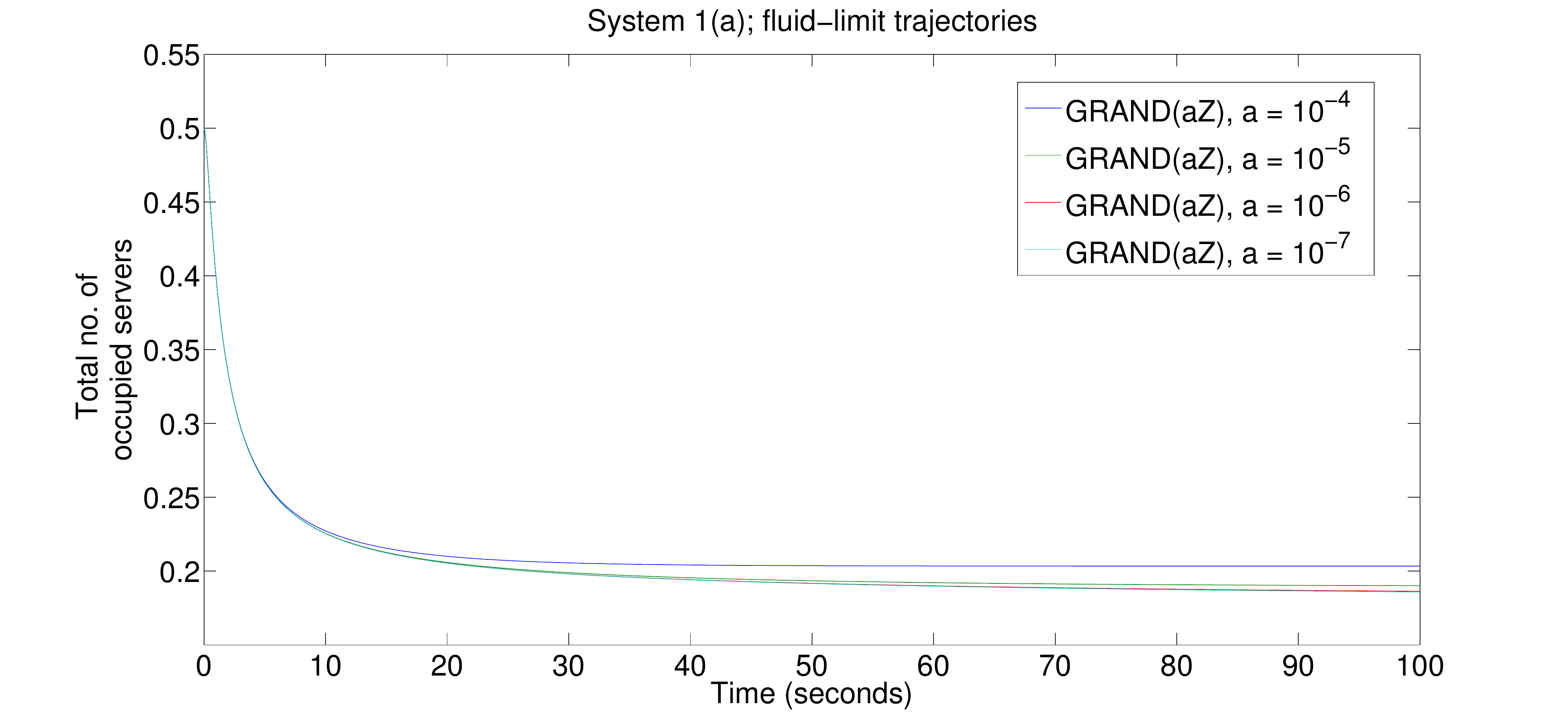}
\caption{Fluid-limit comparison of GRAND($aZ$) for different $a$, in system \ref{fig:vp}. 
$x_{(1, 1)}(0) = 1/2$, and $x_{\bk}(0) = 0$ for $\bk \neq (1, 1)$.}
\label{fig:compare-fluid2}
\end{figure}
\end{center}

The performance benefit of a positive number of zero-servers
seems to be completely absent, even in the fluid limit, in a system with vector packing constraints, 
such as \ref{fig:vp} (see Figure \ref{fig:compare-fluid2}). We suspect that 
this is because in such a system, there does not 
exist ``corner'' configuration such as $(3, 3)$ in \ref{fig:non-vp}; thus, 
``safety stocks'' induced by a positive number of zero-servers 
are not necessary to help the system to ``escape'' from ``bad'' corner configurations.

\paragraph{Summary of observations from simulation and numerical experiments.}
Our discrete-event simulations and the fluid-limit dynamics 
suggest the following. Increasing the number of zero-servers in a GRAND algorithm
makes the convergence to steady state faster, but makes the steady-state performance
(in the sense of the number occupied servers) worse. There are system structures, scales and parameter
ranges for  which increasing the number of zero-servers improves performance during a transition to steady state.
However, it is our observation that for systems of practical scale, a positive number of zero-servers brings no
performance benefit, namely GRAND($0$) performs best, even during transient time intervals.
For systems with the set of feasible configurations defined by vector-packing constraints,
GRAND($0$) appears to perform best for all scales and parameter ranges.

\begin{table}
\caption{Summary of settings} \label{tab:summary}
\begin{center} 
\begin{tabular}{|c|c|c|}
\hline 
  & $r = 1000$ & $r = 10000$ \\
\hline
\ref{fig:vp} & Figure \ref{fig:compare-vp-1000}  & Figure \ref{fig:compare-vp-10000} \\
\hline
\ref{fig:non-vp} & Figure \ref{fig:compare-corner-1000} & Figure \ref{fig:compare-corner-10000} \\
\hline 
\end{tabular} 
\end{center}
\end{table}

\section{Discussion and further work}
\label{sec-further-work}

We introduce the GRAND algorithm, 
and prove that its subclass 
-- GRAND($aZ$) with $a > 0$ -- is asymptotically 
optimal in the sense of minimizing the number 
of occupied servers, as the system scale increases to infinity and
$a \to 0$. 
We have also performed extensive simulations 
to examine various versions of GRAND algorithm, 
and observe a tradeoff between convergence time 
and steady-state performance as we vary the number of zero-servers. 
In principle, we see some performance benefit 
of having a positive number of zero-servers in systems of a very large scale in
transient time periods, 
but in scenarios of practical scale, our simulations 
suggest that GRAND($0$) gives the best performance.
This makes GRAND($0$)  very attractive for practical applications, due to its extreme simplicity.

We now discuss possible extensions of the theoretical results proved in this paper. 
One direction is to analyze versions of GRAND different from GRAND($aZ$).
Main results of this paper concern GRAND($aZ$) with a fixed $a > 0$.
They are fluid-scale results. Namely,
if we consider GRAND($aZ$) even with very small fixed $a > 0$, 
then loss from optimality (in the sense of the number of occupied servers)
is $O(r)$, as $r \to \infty$.  
A natural extension is to consider a GRAND algorithm
with the number of zero-servers being an increasing, but sublinear,
function of $Z$, for example $X_{\bZero} = Z^p$, with $p \in (0, 1)$.
It is natural to expect that this algorithm will have $o(r)$ error, and therefore
we make the following

\begin{conjecture}
\label{conj-sublinear}
Let $p \in (0, 1)$. Consider a sequence of systems operating 
under the GRAND($Z^p$) algorithm, that is 
$X_{\bZero}(t) = \lceil Z(t)^p \rceil$. 
For each $r$, let $\bx^r(\infty)$ denote the random state 
of the fluid-scaled process in the stationary regime. 
Then as $r \rightarrow \infty$, 
\[
d(\bx^r(\infty), \cx^*) \Rightarrow 0.
\]
\end{conjecture}

Our simulations in Section \ref{sec-sim} suggest that the steady-state performance of GRAND
improves as the number of zero-servers decreases. It is then also natural 
to make the following (far reaching) 
conjecture on the asymptotic optimality of GRAND($0$).

\begin{conjecture}
\label{conj}
Consider a sequence of systems operating under the GRAND$(0)$ algorithm. 
For each $r$, let $\bx^r(\infty)$ denote the random state 
of the fluid-scaled process in the stationary regime. Then as $r \rightarrow \infty$, 
\[
d(\bx^r(\infty), \cx^*) \Rightarrow 0.
\]
\end{conjecture}
Proving Conjectures~\ref{conj-sublinear} and \ref{conj} is a potential subject of future work.

Another possible direction is the extension of the GRAND($aZ$) asymptotic optimality --
specifically Theorem \ref{th-grand-fluid} -- to a system with more general service time distributions.
As we explained in Remark 3 in Section \ref{sec-opt-point-interpretation}, 
such a generalization is very natural to expect.

\end{document}